\documentclass[reqno, 12pt]{amsart}
\usepackage{amsfonts}
\usepackage{}
\usepackage{mathrsfs}
\usepackage{amssymb,url}
\usepackage{cite}
\date{}

\allowdisplaybreaks[4]
\theoremstyle{plain}
\newtheorem{thm}{Theorem}[section]

\newtheorem{lem}[thm]{Lemma}

\theoremstyle{remark}

\theoremstyle{definition}

\usepackage[colorlinks,linkcolor=blue,urlcolor=red,linktocpage]{hyperref}

\makeatletter
    
    \newcommand{\Rmnum}[1]{\expandafter\@slowromancap\romannumeral #1@}
  \makeatother

\numberwithin{equation}{section}

\begin{document}
\title
{$L_p$-Blaschke Valuations}

\author[Jin Li]{Jin Li}
\address[Jin Li]{Department of Mathematics, Shanghai University, Shanghai 200444, China}
 \email{\href{mailto: Jin Li
<lijin2955@gmail.com>}{lijin2955@gmail.com}}

\author[Shufeng Yuan]{Shufeng Yuan}
\address[Shufeng Yuan]{Department of Mathematics, Shanghai University, Shanghai 200444, China}
\email{\href{mailto: Shufeng Yuan
<yuanshufeng2003@163.com>}{yuanshufeng2003@163.com}}

\author[Gangsong Leng]{Gangsong Leng}
\address[Gangsong Leng]{Department of Mathematics, Shanghai University, Shanghai 200444, China}
\email{\href{mailto: Gangsong Leng
<gleng@staff.shu.edu.cn>}{gleng@staff.shu.edu.cn}}

\begin{abstract}
In this article, a classification of continuous,
linearly intertwining, symmetric $L_p$-Blaschke ($p>1$) valuations is established as an extension of Haberl's work on Blaschke valuations.
More precisely, we show that for dimensions $n
\geq 3$, the only continuous,
 linearly intertwining, normalized symmetric $L_p$-Blaschke valuation is the normalized $L_p$-curvature image operator, while for dimension $n = 2 $, a rotated normalized $L_p$-curvature image operator is an only additional one. One of the advantages of our approach is that we deal with normalized symmetric $L_p$-Blaschke valuations, which makes it possible to handle the case $p=n$. The cases where $p \not =n$ are also discussed by studying the relations between symmetric $L_p$-Blaschke valuations and normalized ones.
\end{abstract}

\subjclass[2010]{52B45, 52A20}

\keywords{normalized $L_p$-Blaschke valuation, normalized $L_p$-curvature image, $L_p$-Blaschke valuation, $L_p$-curvature image}

\thanks{The authors would like to acknowledge the support
 from the National Natural Science Foundation of China (11271244), Shanghai Leading Academic Discipline Project
(S30104), and Innovation Foundation of Shanghai University (SHUCX120121).}

\maketitle
\section{Introduction}
A \emph{valuation} is a function $Z: \mathcal {Q} \rightarrow \langle \mathcal {G},+ \rangle$ defined on a class of subsets of $\mathbb{R}^n$ with values in an Abelian semigroup $\langle \mathcal {G},+ \rangle$ which satisfies
\begin{align}\label{11}
    Z(K \cup L) + Z(K \cap L) = ZK +ZL,
\end{align}
whenever $K,L,K \cup L,K \cap L \in \mathcal {Q}$.  In recent years, important new results
on the classification of valuations on the space of convex bodies
have been obtained.  The starting point for a systematic investigation of general valuations was
Hadwiger's \cite{Had2} fundamental characterization of the linear combinations of
intrinsic volumes as the continuous valuations that are rigid motion invariant (see \cite{A1,A2,B1,Lud4} for recent important variants).
Its beautiful applications in
integral geometry and geometric probability are described in Hadwiger's book \cite{Had52} and Klain
and Rota's recent book \cite{{Kla97}}.\par

Excellent surveys on the
history of valuations from Dehn's solution of Hilbert's third
problem to approximately 1990 are in McMullen and
Schneider \cite{McM83} or McMullen \cite{McM93}.

First results on convex body valued valuations were obtained by Schneider \cite{Schn1} in the
1970s, where the addition of convex bodies in (\ref{11}) is Minkowski sum.
In recent years, the investigations of convex and star body valued valuations
gained momentum through a series of articles by Ludwig \cite{Lud1,Lud2,ludwigMin,Lud3} (see also \cite{Hab1,Hab2,haberl1,Hab5,Hab3,Par1,Par2,Schn,S1,S2}).
A very recent
development in this area explores the connections between these valuations
and the theory of isoperimetric inequalities (see e.g., \cite{Hab4,S3,Schu}).

Assuming compatibility with the general linear group, Ludwig
\cite{ludwigMin} obtained a complete classification of \emph{$L_p$-Minkowski
valuations}, i.e., valuations where the addition in (\ref{11}) is
$L_p$-Minkowski sum.  Her results establish simple characterizations
of fundamental operators like the projection or centroid body
operator.  \text{Haberl} \cite{haberl1} established a classification of all
continuous symmetric Blaschke valuations, where addition in (\ref{11})
is Blaschke sum ``$\#$", compatible with the
general linear group. For $n \geq 3$, the only two examples of such
valuations are a scalar multiple of the curvature image operator and the Blaschke symmetral $ZK = K \# (-K)$.
For $n=2$, Blaschke sum coincides with Minkowski sum, a classification is provided by \text{Ludwig's} results \cite{ludwigMin}.

In this paper, we extend Haberl's \cite{{haberl1}} results in the
context of the $L_{p}$-Brunn-Minkowski theory when $p>1$ for $n \geq 2$.  To treat the case that $p=n$ when $n$ is not even at the same time as the case for general $p>1$, we deal with normalized symmetric $L_p$-Blaschke valuations (that is the addition in (\ref{11}) is normalized $L_p$-Blaschke sum).  For $n \geq 3$, the only example (up to a dilation) of a continuous, linearly intertwining, normalized symmetric $L_p$-Blaschke valuation is the normalized $L_p$-curvature image operator.  For $n = 2 $, the rotation of the normalized $L_p$-curvature image operator by an angle $\pi /2$ is the only additional example.  As by-products, by the relationship between symmetric $L_p$-Blaschke valuations and corresponding normalized case, we also classify continuous, linearly intertwining, symmetric $L_p$-Blaschke valuations for $p \neq n$.

Since the classification of $L_p$-Blaschke valuations is based on Ludwig's results \cite{ludwigMin}, some other classifications of Minkowski valuations should be remarked here.  Schneider and Schuster \cite{Schn} and Schuster \cite{S1} classified some rotation covariant Minkowski valuations.  Schuster and Wannerer \cite{S2} classified $GL(n)$ contravariant Minkowski valuations without any restrictions on their domain.  Very recently, Haberl \cite{Hab5} showed that the homogeneity assumptions of $p=1$ in Ludwig \cite{ludwigMin} are not necessary, and Parapatits \cite{Par1,Par2} showed that the homogeneity assumptions of $p>1$ in Ludwig \cite{ludwigMin} are also not necessary.  But the homogeneity assumptions are still needed in this paper.

In order to state the main result, we collect some notation. Let $\mathcal{K}^n$ be the space of \emph{convex bodies}, i.e., nonempty, compact, convex subsets of $\mathbb{R}^{n}$, endowed with Hausdorff metric. We denote by
$\mathcal{K}_o ^n$ the set of $n$-dimensional convex bodies which contain the origin, and by
$\overline{\mathcal{K}}_o ^n$ the set of convex bodies which contain
the origin. The set of $n$-dimensional origin-symmetric convex bodies is denoted by $\mathcal{K}_c ^{n}$.

We will always assume that $p \in \mathbb{R}$ and $p > 1$ in this paper, unless noted otherwise.

In \cite{lutwak93}, Lutwak introduced the notion of the $L_p$-surface area measure $S_p(K,\cdot)$ and posed the even
$L_p$-Minkowski problem: given an even Borel measure $\mu$ on the unit sphere $S^{n-1}$,
does there exist an $n$-dimensional convex body $K$ such that $\mu = S_p(K,\cdot)$?
An affirmative answer was given, if $p \neq n$ and if $\mu$ is not concentrated on any great
subsphere.  For $p \neq n$, using the uniqueness of the even $L_p$-Minkowski problem
  on $\mathcal{K}_c ^{n}$, the \emph{$L_p$-Blaschke sum} $K \# _p L \in \mathcal{K}_c ^n$ of
   $K,L \in \mathcal{K}_c ^n$ was defined by
$S_p(K \# _p L,\cdot) = S_p(K,\cdot) + S_p(L,\cdot)$. Thus $\mathcal{K}_c^n$ endowed with $L_p$-Blaschke sum is an Abelian semigroup which we denote by $\langle \mathcal{K}_c ^n,\# _p \rangle$.

The volume-normalized even $L_p$-Minkowski problem, for which the case $p=n$ can be handled as well,
was introduced and solved by Lutwak, Yang and Zhang \cite{LYZ2}.  If $\mu$ is an even Borel measure
on the unit sphere $S^{n-1}$, then there exists a unique $n$-dimensional origin-symmetric convex
body $\widetilde{K}$ such that
\begin{align}\label{lpmp}
\frac{{S_p (\widetilde{K},\cdot)}}{V(\widetilde{K})} = \mu,
\end{align}
if and only if $\mu$ is not concentrated on any great
subsphere, where $V(\widetilde{K})$ is the volume of $\widetilde{K}$.

The volume-normalized even $L_p$-Minkowski problem also suggests the following composition of bodies in $\mathcal{K}_c ^n$.
For $K,L \in \mathcal{K}_c ^n$, we define
the \emph{normalized $L_p$-Blaschke sum} $K \widetilde{\#}_p L \in \mathcal{K}_c ^n$ by
$$\frac{{S_p (K \widetilde{\#}_p L,\cdot)}}{V(K \widetilde{\#}_p L)} = \frac{{S_p (K,\cdot)}}{V(K)} + \frac{{S_p (L,\cdot)}}{V(L)}.$$
Obviously the existence and uniqueness of $K \widetilde{\#}_p L$ are guaranteed by relation (\ref{lpmp}).  Also $\mathcal{K}_c^n$ endowed with the normalized $L_p$-Blaschke sum is an Abelian semigroup which we denote by $\langle \mathcal{K}_c ^n,\widetilde{\#}_p \rangle$.

We call a valuation $Z : \mathcal{K}_o ^n \rightarrow \langle \mathcal{K}_c ^n,\# _p \rangle$ \emph{symmetric $L_p$-Blaschke valuation}, and a valuation $Z : \mathcal{K}_o ^n \rightarrow \langle \mathcal{K}_c ^{n}, \widetilde{\#}_p \rangle$ \emph{normalized symmetric $L_p$-Blaschke valuation}.

A convex body $K$, which contains the origin in its interior, is said to have a \emph{$L_p$-curvature function} $f_p(K,\cdot): S^{n-1} \to \mathbb{R}$, if $S_p(K,\cdot)$ is absolutely continuous with respect to spherical Lebesgue measure $\sigma$, and
\begin{align*}
\frac{dS_p(K,\cdot)}{d \sigma (\cdot)} = f_p(K,\cdot)
\end{align*}
almost everywhere with respect to $\sigma$.

For $p \geq 1$, and $p \neq n$, the \emph{symmetric $L_p$-curvature image} $\Lambda _c ^p K$ of $K \in \mathcal{K}_o ^n$ is defined as the unique body in $\mathcal{K}_c ^n$ such that
$$f_p(\Lambda _c ^p K,\cdot)=\frac{{1}}{2} \rho (K,\cdot)^{n+p} + \frac{{1}}{2} \rho (-K,\cdot)^{n+p},$$
where $\rho _K (\cdot)= \rho (K,\cdot) : S^{n-1} \to \mathbb{R}$ is the radial function of $K$, i.e., $\rho (K,u)= \max \{\lambda >0 : \lambda u \in K \}$.  When $p=1$, this is the classical curvature image operator, a central notion in the affine geometry of convex bodies;
see e.g., \cite{Lei1,Li1,Lut3,Lut1,Lut2,Lut4}.  When $p>1$, it should be noticed that the
definition of the $L_p$-curvature image operator here differs from the definition of the Lutwak \cite{Lut5}.

For $p \geq 1$, the \emph{normalized symmetric $L_p$-curvature image} $\widetilde{\Lambda} _c ^p K$ of $K \in \mathcal{K}_o ^n$ is defined as the unique body in $\mathcal{K}_c ^n$ such that
$$\frac{{f_p(\widetilde{\Lambda} _c ^p K,\cdot)}}{V(\widetilde{\Lambda} _c ^p K)}=(\frac{{1}}{2} \rho (K,\cdot)^{n+p} + \frac{{1}}{2} \rho (-K,\cdot)^{n+p}).$$

\textbf{Remark:} By the uniqueness of the even $L_p$-Minkowski problem and the volume-normalized even $L_p$-Minkowski problem, if $p \geq 1$, and $p \neq n$, it follows that
$$V(\widetilde{\Lambda} _c ^pK)^{1/(p-n)} \widetilde{\Lambda} _c ^p K = \Lambda _c ^p K.$$

An operator $Z : \mathcal {Q} \to \langle \mathfrak{P}(\mathbb{R}^n),+ \rangle$, where $\mathfrak{P}(\mathbb{R}^n)$ denotes the power set of $\mathbb{R}^n$, is called \emph{$SL(n)$ covariant}, if
\begin{align*}
Z(\phi K) = \phi ZK
\end{align*}
for every $K \in \mathcal {Q}$ and $\phi \in SL(n)$.  It is called \emph{$SL(n)$ contravariant}, if
\begin{align*}
Z(\phi K) = \phi ^{-t} ZK
\end{align*}
for every $K \in \mathcal {Q}$ and $\phi \in SL(n)$.  Here, $\phi ^{-t}$ denotes the inverse of the transpose of $\phi$.  We call Z \emph{homogeneous of degree $q \in \mathbb{R}$}, if
\begin{align*}
Z(\lambda K) = \lambda ^q ZK
\end{align*}
for every $K \in \mathcal {Q}$ and $\lambda >0$, and we call $Z$ homogeneous if it is homogeneous for some $q \in \mathbb{R}$.
If $Z$ is homogeneous and $SL(n)$ covariant or contravariant, then we call it \emph{linearly intertwining}.

Our main results are the following two theorems.
\begin{thm}\label{1.1}
    Let $n \geq 2$. For $p > 1$ and $p$ not an even integer, the operator $Z : \mathcal{K}_o ^n \rightarrow \langle \mathcal{K}_c ^n,\widetilde{\#}_p \rangle$ is a continuous, homogeneous, $SL(n)$ contravariant valuation, if and only if there exists a constant $c>0$ such that
$$ZK=c \widetilde{\Lambda} _c ^p K$$
for every $K \in \mathcal{K}_o ^n$.
\end{thm}

\begin{thm}\label{1.2}
Let $n \geq 3$. For $p >1$ and $p$ not an even integer, there are no continuous, homogeneous, $SL(n)$ covariant normalized symmetric $L_p$-Blaschke valuations on $\mathcal{K}_o ^n$.

For $p >1$ and $p$ not an even integer, the operator $ Z : \mathcal{K}_o ^2 \rightarrow \langle \mathcal{K}_c ^2,\widetilde{\#}_p \rangle$ is a continuous, homogeneous, $SL(2)$ covariant valuation, if and only if there exists a constant $c>0$ such that
$$ZK=c \psi _{\pi /2} \widetilde{\Lambda} _c ^p K$$
for every $K \in \mathcal{K}_o ^2$.  Here $\psi _{\pi /2}$ is the rotation by an angle $\pi /2$.
\end{thm}

Theorem 1.1 and 1.2 establish a classification of continuous, linearly intertwining, normalized symmetric $L_p$-Blaschke valuations on $\mathcal{K}_o ^n$ when $p>1$ and $p$ is not an even integer. For $p=1$, Haberl \cite{{haberl1}} obtained a complete classification of continuous, linearly intertwining symmetric Blaschke valuations and we can easily get the corresponding results in the normalized case by reversing the process of Theorem \ref{th5.3} and Theorem \ref{th5.5}. Therefore we state the results here only for $p > 1$.

In Section 2, some preliminaries are given. The aim of Section 3 is to derive the
characterizing properties (stated in Theorem 1.1) of the normalized symmetric
$L_p$-curvature image operator $\widetilde{\Lambda} _c ^p$. In Section 4, Lemma 4.1 - Lemma 4.5 generate a homogeneous, $SL(n)$ covariant $L_p$-Minkowski valuation on $\overline{\mathcal{K}}_o ^n$ by a continuous, homogeneous, $SL(n)$ contravariant normalized symmetric $L_p$-Blaschke valuation on $\mathcal{K}_o ^n$.  Using properties of the support set of the $L_p$-projection bodies established in Lemma \ref{pgz} and characterization theorems of $L_p$-Minkowski valuations \cite{ludwigMin}, we classify continuous, homogeneous, $SL(n)$ contravariant normalized symmetric $L_p$-Blaschke valuations. In a similar way, we also classify continuous, homogeneous, $SL(n)$ covariant normalized symmetric $L_p$-Blaschke valuations.
In Section 5, from the relationship between normalized symmetric $L_p$-Blaschke valuations and symmetric $L_p$-Blaschke valuations (Lemma \ref{th5.1} and Lemma \ref{th5.2}), we also classify continuous, linearly intertwining, symmetric $L_p$-Blaschke valuations on $\mathcal{K}_o ^n$  for $p \neq n$ (see Theorem \ref{th5.3} and Theorem \ref{th5.5}).
\section{Preliminaries}
We work in Euclidean $n$-space $\mathbb{R}^n$ with $n \geq 2$. Let $\{e_i\}, i=1,\cdots, n$ be the standard basis of $\mathbb{R}^n$. The usual scalar product of two vectors $x$ and $y \in \mathbb{R}^n$
shall be denoted by $x \cdot y$. For $u \in S^{n-1}$, $u^- = \{x \in \mathbb{R}^n : x \cdot u \leq 0\}$,
$u^+ = \{x \in \mathbb{R}^n : x \cdot u \geq 0\}$ and $u^ \perp  = \{x \in \mathbb{R}^n : x \cdot u = 0\}$.   The
convex hull of a set $A \subset \mathbb{R}^n$ will be denoted by $[A]$. To shorten the notation we write
$[A,\pm x_1,\cdots,\pm x_m]$ instead of $[A \cup \{x_1,-x_1,\cdots,x_m,-x_m\}]$ for $ A \subset \mathbb{R}^n$, $m \in \mathbb{N}$, and
$x_1,\cdots, x_m \in \mathbb{R}^n$.  In $\mathbb{R}^2$, we write $\psi _{\pi/2}$ for the rotation by an angle $\pi /2$.

The \emph{Hausdorff distance} of two convex bodies $K,L$ is defined as
$d(K,L)= \mathop {\max }\limits_{u \in S^{n - 1} }$ $ |h_K(u) -
h_L(u)|$, where $h_K : \mathbb{R}^n \to \mathbb{R}$ is the \emph{support function} of
$K \in \mathcal{K}^n$, i.e., $h_K(x)=\max \{x \cdot y : y \in K \}$.
Sometimes we also write $h_K(\cdot)$ as $h(K,\cdot)$. If $f:\mathbb{R}^n \to \mathbb{R}$ is a sublinear function (i.e., $f(\lambda x) = \lambda f(x)$ for every $\lambda \geq 0$ and $x \in \mathbb{R}^n$; $f(x+y) \leq f(x)+f(y)$ for every $x,y \in \mathbb{R}^n$), then there exists a unique convex body $K$ such that $f=h_K$.

Let $S(K,\cdot)$ be the classical surface area measure of a convex body $K$. If $K$ contains the origin in its interior, the Borel measure $S_p(K,\cdot) = h_K(\cdot)^{1-p}S(K,\cdot)$ on $S^{n-1}$ is the \emph{$L_p$-surface area measure} of $K$.

For $K,L \in \mathcal{K}^n$ and $\alpha , \beta \geq 0$ (not both $0$), the \emph{Minkowski linear combination} $\alpha K+ \beta L$ is defined by $\alpha K+ \beta L=\{\alpha x+ \beta y:x \in K,y \in L\}$. For $K,L \in \overline{\mathcal{K}}_o ^n$ and $\alpha , \beta \geq 0$, the \emph{$L_p$-Minkowski linear combination} $\alpha \cdot K+_p \beta \cdot L$ (not both $0$) is defined by $h(\alpha \cdot K+_p \beta \cdot L,u)^p = \alpha h(K,u)^p+ \beta h(L,u)^p$ for every $u \in S^{n-1}$. Note that "$\cdot$" rather than "$\cdot _p$" is written for $L_p$-Minkowski scalar multiplication. This should create no confusion. Also note that the relationship between $L_p$-Minkowski and Minkowski scalar multiplication is $\alpha \cdot K = \alpha ^{1/p} K$.

For $p \geq 1$, the $L_p$-mixed volume $V_p(K,L)$ of the convex bodies $K,L$ containing the origin in their interiors was defined in \cite{lutwak93} by
\begin{align*}
\frac{n}{p}V_p(K,L)= \mathop {\lim }\limits_{\varepsilon \to 0^ + }\frac{V(K+_p \varepsilon \cdot L)-V(K)}{\varepsilon},
\end{align*}
where the existence of this limit was demonstrated in \cite{lutwak93}. Obviously, for each $K$, $V_p(K,K) = V(K)$. It was also shown in \cite{lutwak93} that the $L_p$-mixed volume $V_p$ has the following integral representation:
\begin{align*}
V_p(K,L)=\frac{1}{n} \int _{S^{n-1}} h(L,u)^p dS_p(K,u).
\end{align*}

For $p\geq 1$, the \emph{$L_p$-cosine transform} of a finite, signed Borel measure $\mu$ on $S^{n-1}$ is defined by
$$C_p \mu (x) = \int _{S^{n-1}} | x \cdot v |^p d \mu (v),~~~~x \in \mathbb{R}^n.$$
Similarly, the \emph{$L_p$-cosine transform} of a Borel measurable function $f$ on $S^{n-1}$ is defined by
$$(C_p f) (x) = \int _{S^{n-1}} |x \cdot v|^p f(v) d \sigma (v),~~~~x \in \mathbb{R}^n,$$
where $\sigma $ is the spherical Lebesgue measure.
An important property of this integral transform is the following injectivity behavior.  If $p$ is not an even integer, and $\mu$ is a signed finite even Borel measure, then
\begin{align}\label{inj}
    \int _{S^{n-1}} |u \cdot v|^p d \mu (v) = 0 ~\text{for all}~ u \in S^{n-1} \Rightarrow \mu =0.
\end{align}
(see e.g.,  Koldobsky \cite{K1,K2}, Lonke \cite{L1}, Neyman
\cite{N1}, and Rubin \cite{R1,Rubin1}.)

For $p \geq 1$, the \emph{$L_p$-projection body}, $\Pi _p K$, of a convex body $K$ containing the origin in its interior is the origin-symmetric convex body whose support function is defined by
\begin{align*}
h(\Pi _p K,u)^p = \int _{S^{n-1}} |u \cdot v|^p dS_p(K,v)
\end{align*}
for every $u \in S^{n-1}$.  The notion of the $L_p$-projection body (with a different normalization)
was introduced by Lutwak, Yang and Zhang \cite{LYZ1}.

It is proved in \cite{LYZ1} that
\begin{align*}
\Pi _p \phi K = |\det \phi|^ {1/p} \phi ^{-t} \Pi _p K
\end{align*}
for every $\phi \in GL(n)$. Then we immediately get
\begin{align}\label{Cp}
    C_p S_p(\phi K, \cdot)(x)=|\det \phi|C_p S_p(K,\cdot)(\phi ^ {-1} x),
\end{align}
and
\begin{align}\label{ncp}
    C_p \frac{S_p(\phi K, \cdot)}{V(\phi K)}(x)=C_p \frac{S_p(K,\cdot)}{V(K)}(\phi ^ {-1} x).
\end{align}

The notion of the \emph{$L_p$-centroid body} was introduced by Lutwak, Yang
and Zhang \cite{LYZ1}: For each compact star-shaped (about the
origin) $K$ in $\mathbb{R}^{n}$ and for $ p\geq1$, the $L_{p}$-centroid
body $\Gamma _p K$ is defined by
\begin{align}\label{CCp}
h(\Gamma _p K,u) = (\frac{{1}}{c_{n,p} V(K)} \int _K |x \cdot u|^p dx)^{1/p}
\end{align}
for every $u \in S^{n-1}$, where the constant $c_{n,p}$ is chosen so that $\Gamma _p B = B$.  For $p = 2$, the $\Gamma _2$-centroid body is the Legendre ellipsoid of classical mechanics.  It is easy to see that
\begin{align}\label{2.10}
\Gamma _p \phi K = \phi \Gamma _p K
\end{align}
for every $\phi \in GL(n)$.
We also can rewrite relation (\ref{CCp}) for the $L_p$-cosine transform:
\begin{align}\label{2.11}
 h(\Gamma _p K,u)^p
 &= \frac{{1}}{(n+p) c_{n,p} V(K)} (C_p \rho _K^{n+p}) (u)\nonumber \\
&= \frac{{1}}{(n+p) c_{n,p} V(K)} (C_p (\frac{{1}}{2}\rho _K^{n+p} +
\frac{{1}}{2}\rho _{-K}^{n+p})) (u).
\end{align}

\section{Normalized symmetric $L_p$-curvature images}
In this section, we will show that the normalized symmetric $L_p$-curvature image operator $\widetilde{\Lambda} _c ^p$ is a continuous, homogeneous, $SL(n)$ contravariant normalized symmetric $L_p$-Blaschke valuation.

We remark that a valuation $Z : \mathcal {Q} \to \langle \mathfrak{P}(\mathbb{R}^n),+ \rangle$ is $SL(n)$ covariant and homogeneous of degree $q$ if and only if it satisfies
\begin{align}\label{coho}
    Z(\phi K) = (\det \phi)^ \frac{{q - 1}}{n} \phi ZK
\end{align}
for every $K \in \mathcal {Q}$ and $\phi \in GL(n)$ with positive determinant.  Similarly, a valuation $Z$ is $SL(n)$
contravariant and homogeneous of degree $q$ if and only if it satisfies
\begin{align}\label{ctho}
    Z(\phi K) = (\det \phi)^ \frac{{q + 1}}{n} \phi ^{-t} ZK
\end{align}
for every $K \in \mathcal {Q}$ and $\phi \in GL(n)$ with positive determinant.

To prove that $\widetilde{\Lambda} _c ^p$ is a continuous valuation, we will firstly show the following lemma.
\begin{lem}\label{th4.2}
If $K_i,K \in \mathcal{K}_c ^n,~i=1,2,\cdots$, such that $\frac{{S_p(K_i,\cdot)}}{V(K_i)} \to \frac{{S_p(K,\cdot)}}{V(K)}$ weakly, then $K_i \to K$.
\begin{proof}
Firstly, we want to show that $\{K_i\}$ has a subsequence, $\{ K_ {i_j} \}$, converging to an origin-symmetric convex body containing the origin in its interior (the proof is similar as \cite[Theorem 2]{LYZ2}).

Define $f_K(u)$ by
$$f_K(u) ^p= \frac{{1}}{n} \int _{S^{n-1}} |u \cdot v|^p \frac{{dS_p(K,v)}}{V(K)}.$$
Thus $f_K(u)$ is a support function of some convex body.
Since $\frac{{S_p(K,\cdot)}}{V(K)}$ is not concentrated on any great subsphere, $f_K(u)>0$ for every $u \in S^{n-1}$.  By the continuity of $f_K(u)$ on the compact set $S^{n-1}$, there exist two constants $a,b>0$, such that $\frac{{1}}{2}a \geq f_K(u) \geq 2b$ for every $u \in S^{n-1}$.  Since $\frac{{S_p(K_i,\cdot)}}{V(K_i)} \to \frac{{S_p(K,\cdot)}}{V(K)}$ weakly, we get $f_{K_i}(u) \to f_K(u)$. The convergence is uniform in $u \in S^{n-1}$ by \cite[Theorem 1.8.12]{Schn2}. Hence $a \geq f_{K_i} \geq b$ for sufficiently large $i$ uniformly.

In order to show $K_i$ is uniformly bounded, define real numbers $M_i$, and vectors $u_i \in S^{n-1}$ by
\begin{align*}
M_i = \mathop {\max }\limits_{u \in S^{n-1}} h(K_i,u) = h(K_i,u_i).
\end{align*}
Now, $M_i [-u_i,u_i] \subset K_i$.  Hence $M_i |u_i \cdot v| \leq h(K_i,v)$ for every $v \in S^{n-1}$.  Thus,
\begin{align*}
M_i ^p b^p & \leq M_i ^p \frac{{1}}{n} \int _{S^{n-1}} |u_i \cdot v|^p \frac{{dS_p(K_i,v)}}{V(K_i)} \\
&\leq \frac{{1}}{n} \int _{S^{n-1}} h(K_i,v)^p \frac{{dS_p(K_i,v)}}{V(K_i)}
= \frac{{V_p( K_i,K_i)}}{V(K_i)}
= 1
\end{align*}
for sufficiently large $i$. Hence $K_i$ is uniformly bounded. By the Blaschke selection theorem, there exists a subsequence $\{K_ {i_{j}} \}$ converging to a convex body, say $K'$. Since $K_ {i_{j}}$ are origin-symmetric, $K'$ is origin-symmetric. Define real numbers $m_i$, and vectors $u_i ' \in S^{n-1}$ by
\begin{align*}
m_i = \mathop {\min }\limits_{u \in S^{n-1}} h(K_i,u) = h(K_i,u_i '),
\end{align*}
The property $a \geq f_{K_i}$ for sufficiently large $i$ uniformly, together with Jensen's inequality, shows that \begin{align*}
a & \geq (\frac{{1}}{n} \int _{S^{n-1}} |u_i ' \cdot v|^p \frac{{dS_p(K_i,v)}}{V(K_i)})^\frac{{1}}{p}
= (\frac{{1}}{n} \int _{S^{n-1}} (\frac{{|u_i ' \cdot v|}}{h(K_i,v)}) ^p \frac{{h(K_i,v)dS(K_i,v)}}{V(K_i)})^\frac{{1}}{p} \\
& \geq \frac{{1}}{n} \int _{S^{n-1}} \frac{{|u_i ' \cdot v|}}{h(K_i,v)} \frac{{h(K_i,v)dS(K_i,v)}}{V(K_i)}
= \frac{{2V(K_i | (u_i ') ^ \perp)}}{nV(K_i)}.
\end{align*}
Since $K_i$ is contained in the right cylinder $K_i | (u_i ') ^ \perp \times m_i [-u_i ',u_i ']$, we have \\
$2 m_i V(K_i | (u_i ') ^ \perp) \geq V(K_i)$.  Thus,
\begin{align*}
a \geq \frac{{2V(K_i | (u_i ') ^ \perp)}}{nV(K_i)} \geq \frac{{1}}{nm_i},
\end{align*}
which shows $m_i \geq \frac{{1}}{na}$for sufficient large $i$.  Hence
\begin{align*}
\frac{{1}}{na} B \subseteq K',
\end{align*}
where $B$ is the unit ball in $\mathbb{R}^n$.  Thus, $K'$ contains the origin in its interior. The first step is complete.

Next, we argue the assertion by contradiction.  Assume $K_i \nrightarrow K$, then there exists a subsequence, $\{K_{i_j} \}$, such that $d(K_{i_j},K) \geq \varepsilon$ for a suitable $\varepsilon >0$.  Since $\{K_{i_j} \}$ also satisfies the condition of this Lemma, from the conclusion above, there exists a subsequence of $\{ K_ {i_j} \}$, say $\{K_ {i_{j_k}} \}$, converging to an origin-symmetric convex body, say $K'$, containing the origin in its interior.  Thus,
$\frac{{S _{p} (K_ {i_{j_k}}, \cdot)}}{V(K_{i_{j_k}})} \to \frac{{S _{p} (K', \cdot)}}{V(K')}$ weakly.  By the uniqueness of weak convergence and the normalized even $L_p$-Minkowski problem, we get $K_{i_{j_k}} \to K' = K$.  That is a contradiction.
\end{proof}
\end{lem}

\begin{thm}\label{th3.1}
The normalized symmetric $L_p$-curvature image operator $\widetilde{\Lambda} _c ^p:\mathcal{K}_o ^n \to \langle \mathcal{K}_c ^n, \widetilde{\#}_p \rangle$ is a continuous, $SL(n)$ contravariant valuation which is homogeneous of degree $-\frac{n}{p}-1$.  Moreover, $\psi _{\pi /2} \widetilde{\Lambda} _c ^p:\mathcal{K}_o ^2 \to \langle \mathcal{K}_c ^2,\widetilde{\#}_p \rangle$ is a continuous, $SL(2)$ covariant valuation which is homogeneous of degree $-\frac{2}{p}-1$.
\begin{proof}
To prove that $\widetilde{\Lambda} _c ^p$ is a normalized symmetric $L_p$-Blaschke valuation, we just need to show
\begin{align}\label{41}
\frac{{S_p(\widetilde{\Lambda} _c ^p (K \cup L),\cdot)}}{V(\widetilde{\Lambda} _c ^p(K \cup L))}
 + \frac{{S_p(\widetilde{\Lambda} _c ^p(K \cap L),\cdot)}}{V(\widetilde{\Lambda} _c ^p(K \cap L))}
 = \frac{{S_p(\widetilde{\Lambda} _c ^p K,\cdot)}}{V(\widetilde{\Lambda} _c ^p K)}
 + \frac{{S_p(\widetilde{\Lambda} _c ^p L,\cdot)}}{V(\widetilde{\Lambda} _c ^p L)}
\end{align}
for every $K,L,K \cup L,K \cap L \in \mathcal{K}_o ^n$.  Since
\begin{align*}
 \rho (K \cup L,\cdot)^{n+p} + \rho (K \cap L,\cdot)^{n+p} &= \rho (K,\cdot)^{n+p} + \rho (L,\cdot)^{n+p}, \\
 \rho (-(K \cup L),\cdot)^{n+p} + \rho (-(K \cap L),\cdot)^{n+p} &= \rho (-K,\cdot)^{n+p} + \rho (-L,\cdot)^{n+p}
\end{align*}
for every $K,L,K \cup L,K \cap L \in \mathcal{K}_o ^n$, it follows from the definition of $\widetilde{\Lambda} _c ^p$, that the relation (\ref{41}) is true.  Hence the valuation property is established.

To prove homogeneity and $SL(n)$ contravariance of $\widetilde{\Lambda} _c ^p$, by relation (\ref{ctho}), we need to show
\begin{align}\label{42}
\widetilde{\Lambda} _c ^p \phi K = (\det \phi)^ {-1/p} \phi ^{-t} \widetilde{\Lambda} _c ^p K
\end{align}
for every $\phi \in GL(n)$ with positive determinant.  Indeed, the definition of $\widetilde{\Lambda} _c ^p$, the relation (\ref{2.10}), (\ref{2.11}) together with (\ref{ncp}) imply that
\begin{align*}
C_p \frac{{S_p(\widetilde{\Lambda} _c ^p \phi K,\cdot)}}{V(\widetilde{\Lambda} _c ^p \phi K)}(u)
&= (C_p (\frac{{1}}{2} \rho _ {\phi K} ^{n+p} + \frac{{1}}{2} \rho _{-\phi K} ^{n+p}))(u) \\
&= (n+p) c_{n,p} V(\phi K) h(\Gamma _p \phi K,u)^p \\
&= |\det \phi|(n+p) c_{n,p} V(K) h(\Gamma _p K,\phi ^t u)^p \\
&= |\det \phi| C_p \frac{{S_p(\widetilde{\Lambda} _c ^p K,\cdot)}}{V(\widetilde{\Lambda} _c ^p K)}(\phi ^t u) \\
&= C_p \frac{{S_p(|\det \phi| ^{-1/p} \phi ^{-t} \widetilde{\Lambda} _c ^p K,\cdot)}}{V(|\det \phi| ^{-1/p} \phi ^{-t} \widetilde{\Lambda} _c ^p K)}(u).
\end{align*}
The injectivity property (\ref{inj}) and the uniqueness of the volume-normalized even $L_p$-Minkowski problem now imply relation (\ref{42}).

If $K_i \to K$, then $\rho (K_i,\cdot) \to \rho (K,\cdot)$ almost everywhere with respect to spherical Lebesgue measure (see \cite[Lemma 1]{haberl1}).  Hence
$(\frac{{1}}{2} \rho (K_i,\cdot)^{n+p} + \frac{{1}}{2} \rho (-K_i,\cdot)^{n+p}) \to (\frac{{1}}{2} \rho (K,\cdot)^{n+p} + \frac{{1}}{2} \rho (-K,\cdot)^{n+p})$ almost everywhere.  Since $(\frac{{1}}{2} \rho (K_i,\cdot)^{n+p} + \frac{{1}}{2} \rho (-K_i,\cdot)^{n+p})$ are uniformly bounded, $\frac{{S_p(\widetilde{\Lambda} _c ^pK_i,\cdot)}}{V(\widetilde{\Lambda} _c ^pK_i)} \to \frac{{S_p(\widetilde{\Lambda} _c ^pK,\cdot)}}{V(\widetilde{\Lambda} _c ^pK)}$ weakly.  Hence, by Lemma \ref{th4.2}, we get $\widetilde{\Lambda} _c ^pK_i \to \widetilde{\Lambda} _c ^pK$.  Thus, $\widetilde{\Lambda} _c ^p K$ is a continuous valuation.

If $\phi \in SL(2)$, we have $\psi _{\pi /2} \phi ^{-t} \psi _{-\pi /2} = \phi$. Then we get
$$\psi _{\pi /2} \widetilde{\Lambda} _c ^p \phi K = \psi _{\pi /2} \phi ^{-t} \widetilde{\Lambda} _c ^p K
= \psi _{\pi /2} \phi ^{-t} \psi _{-\pi /2} \psi _{\pi /2} \widetilde{\Lambda} _c ^p K = \phi \psi _{\pi /2} \widetilde{\Lambda} _c ^p K$$
for every $K \in \mathcal{K}_o ^n$.
Since the operator $\psi _{\pi /2}$ is continuous, we obtain that $\psi _{\pi /2} \widetilde{\Lambda} _c ^p$ is continuous. Moreover, it is easy to verify that $\psi _{\pi /2} \widetilde{\Lambda} _c ^p$ is a normalized symmetric $L_p$-Blaschke valuation which is homogeneous of degree $-\frac{2}{p}-1$.  Hence, $\psi _{\pi /2} \widetilde{\Lambda} _c ^p$ is a continuous, $SL(2)$ covariant normalized symmetric $L_p$-Blaschke valuation which is homogeneous of degree $-\frac{2}{p}-1$.
\end{proof}
\end{thm}

\section{ Normalized $L_p$-Blaschke Valuations }
In this section, for the contravariant and covariant case respectively, we establish our classification results for continuous, linearly intertwining, normalized symmetric
$L_p$-Blaschke valuations.

We remark firstly the fact that the $SL(n)$ covariance (or contravariance) and homogeneity of a valuation $Z: \overline{\mathcal{K}}_o ^n \rightarrow \langle \mathfrak{P}(\mathbb{R}^n),+ \rangle$ are completely determined by the restriction of $Z$ to $n$-dimensional convex bodies if the Abelian semigroup $\langle \mathfrak{P}(\mathbb{R}^n),+ \rangle$ has the cancellation property.  (Actually this property is generalized from Lemma 4 and Lemma 9 of Haberl \cite{haberl1}, and the proof of this property is almost the same as Haberl's).

\begin{lem}\label{2th1}
    If $Z: \overline{\mathcal{K}}_o ^n \rightarrow \langle \mathfrak{P}(\mathbb{R}^n),+ \rangle$ is a valuation which is $SL(n)$ covariant (or contravariant) and homogeneous
    of degree q on $n$-dimensional convex bodies, and $\langle \mathfrak{P}(\mathbb{R}^n),+ \rangle$ has the cancellation property, then $Z$ is $SL(n)$ covariant (or contravariant respectively) and homogeneous of degree q on $\overline{\mathcal{K}}_o ^n$.
\begin{proof}
    In the covariant case, we have to show
\begin{align}\label{f1}
    Z\phi K=(\det \phi)^ \frac{{q - 1}}{n} \phi ZK
\end{align}
    for every $K \in \overline{\mathcal{K}}_o ^n$ and $\phi \in GL(n)$ with positive determinant.  Let $\dim K = n-k$, where
    $0 \leq k \leq n$.  We prove our assertion by induction on $k$.  Indeed, (\ref{f1}) is true for $k=0$ by assumption.  Assume that
    (\ref{f1}) holds for $(n-k)$-dimensional convex bodies and $\dim K =n-(k+1)$.  Choose $u \notin \text{lin}~K$, where $\text{lin}~K$ denotes the linear hull of $K$.  Clearly $[K,u],~[K,-u],~[K,u,-u]$, $\phi [K,u],~\phi [K,-u],~\phi [K,u,-u]$ are of dimension $n-k$, and
\begin{align*}
    [K,u] \cup [K,-u] = [K,u,-u], &~[K,u] \cap [K,-u] = K,  \\
    \phi [K,u] \cup \phi [K,-u] = \phi [K,u,-u], & ~\phi [K,u] \cap \phi [K,-u] = \phi K.
\end{align*}
    Since Z is a valuation,
\begin{align*}
    Z \phi K + Z \phi [K,u,-u] = Z \phi [K,u] + Z \phi [K,-u].
\end{align*}
    With the induction assumption, we get
\begin{align*}
    Z \phi K + (\det \phi)^ \frac{{q - 1}}{n} \phi Z [K,u,-u]
    = (\det \phi)^ \frac{{q - 1}}{n} \phi Z [K,u] + (\det \phi)^ \frac{{q - 1}}{n} \phi Z [K,-u].
\end{align*}
    So,
\begin{align*}
    (\det \phi)^ {-\frac{{q - 1}}{n}} \phi ^{-1}Z \phi K + Z [K,u,-u]
    = Z [K,u] + Z [K,-u].
\end{align*}
    By the cancellation property of $\langle \mathfrak{P}(\mathbb{R}^n),+ \rangle$, combined with the relation
\begin{align*}
    ZK + Z[K,u,-u] = Z[K,u] + Z [K,-u],
\end{align*}
    we have
\begin{align}
    (\det \phi)^ {-\frac{{q - 1}}{n}} \phi ^{-1}Z \phi K = ZK.
\end{align}
    This immediately proves that (\ref{f1}) holds for bodies of dimension $n-k-1$.

    The contravariant case is proved similarly to the covariant case.
\end{proof}
\end{lem}

Since $\mathcal{K}^n$ endowed with $L_p$-Minkowski sum is an Abelian semigroup which has the cancellation property, we immediately
get the following.
\begin{lem}\label{Min}
    If $Z: \overline{\mathcal{K}}_o ^n \rightarrow \langle \mathfrak{P}(\mathbb{R}^n),+_p \rangle$ is a $L_p$-Minkowski valuation which is $SL(n)$ covariant (or contravariant) and homogeneous of degree q on $n$-dimensional convex bodies, then $Z$ is $SL(n)$ covariant (or contravariant respectively) and homogeneous of degree q on $\overline{\mathcal{K}}_o ^n$.
\end{lem}

\subsection{ The contravariant case }
Firstly, we reduce the possible degrees of homogeneity of continuous, $SL(n)$ contravariant normalized symmetric $L_p$-Blaschke valuations.
\begin{lem}\label{3.1}
If $ Z : \mathcal{K}_o ^n \rightarrow \langle \mathcal{K}_c ^n,\widetilde{\#}_p \rangle$ is a continuous, $SL(n)$ contravariant valuation which is homogeneous of degree $q$, then $q \leq -1$.
\begin{proof}
Suppose $K \in \mathcal{K}_o ^n$ is an arbitrary convex body that $K \cap e_n ^+$ and $K \cap e_n ^-$
are $n$-dimensional.  For every positive $s$ we have
$$[K \cap e_n ^+ ,\pm se_n] \cup [K \cap e_n ^- ,\pm se_n] = [K,\pm se_n],$$
$$[K \cap e_n ^+ ,\pm se_n] \cap [K \cap e_n ^- ,\pm se_n] = [K \cap e_n ^ \perp ,\pm se_n].$$
Since $ Z$ is a normalized symmetric $L_p$-Blaschke valuation, we have
\begin{align}\label{ct1}
    &C_p \frac{{S_p(Z[K \cap e_n ^ \perp ,\pm se_n],\cdot)}}{V(Z[K \cap e_n ^ \perp ,\pm se_n])}(e_1) \nonumber \\
    =& C_p \frac{{S_p(Z[K \cap e_n ^+ ,\pm se_n],\cdot)}}{V(Z[K \cap e_n ^+ ,\pm se_n])}(e_1) + C_p \frac{{S_p(Z[K \cap e_n ^- ,\pm se_n],\cdot)}}{V(Z[K \cap e_n ^- ,\pm se_n])}(e_1) \nonumber \\
    & - C_p \frac{{S_p(Z[K,\pm se_n],\cdot)}}{V(Z[K,\pm se_n])}(e_1).
\end{align}
Define a linear map $\phi$ by
$$\phi e_i = e_i, i=1,\cdots,n-1, \phi e_n =se_n.$$
From the $SL(n)$ contravariance and homogeneity of $Z$ as well as relation (\ref{ctho}) and (\ref{ncp}), we get
\begin{align*}
    C_p \frac{{S_p(Z[K \cap e_n ^ \perp ,\pm se_n],\cdot)}}{V(Z[K \cap e_n ^ \perp ,\pm se_n])}(e_1)
    & = C_p \frac{{S_p(s^{ \frac{{q+1}}{n} } \phi ^{-t}Z[K \cap e_n ^ \perp ,\pm e_n],\cdot)}}{V(s^{ \frac{{q+1}}{n} } \phi ^{-t}Z[K \cap e_n ^ \perp ,\pm e_n])}(e_1)
    \\ & =s^ \frac{{-(q+1)p}}{n} C_p \frac{{S_p(Z[K \cap e_n ^ \perp ,\pm e_n],\cdot)}}{V(Z[K \cap e_n ^ \perp ,\pm e_n])}(\phi ^t e_1).
\end{align*}
Since $|e_1 \cdot u|>0$ for all $u \in S^{n-1} \setminus {e_1 ^ \perp}$, and the $L_p$-surface area measure of $n$-dimensional bodies is not concentrated on any great sphere,
we conclude that
\begin{align*}
   & C_p \frac{{S_p(Z[K \cap e_n ^ \perp ,\pm e_n],\cdot)}}{V(Z[K \cap e_n ^ \perp ,\pm e_n])}(\phi ^t e_1) \\
&= \frac{{1}}{V(Z[K \cap e_n ^ \perp ,\pm e_n])} \int _{S^{n-1}} |e_1 \cdot u|^p dS_p(Z[K \cap e_n ^ \perp ,\pm e_n],u)>0.
\end{align*}
Moreover we have
\begin{align*}
    \begin{array}{*{20}c}
    \hfill \mathop {\lim }\limits_{s \to 0^ + } [K \cap e_n ^+ ,\pm se_n] & = & K \cap e_n ^+, \hfill \\
    \hfill \mathop {\lim }\limits_{s \to 0^ + } [K \cap e_n ^- ,\pm se_n] & = & K \cap e_n ^-, \hfill \\
    \hfill \mathop {\lim }\limits_{s \to 0^ + } [K , \pm se_n] & = & K.  \hfill
    \end{array}
\end{align*}
Hence the continuity of $Z$ and volume, together with the weak continuity of $L_p$-surface area measures imply that the right side of (\ref{ct1})
converges to a finite number as $s \to 0^+$.  This implies $\frac{{-(q+1)p}}{n} \geq 0$, so $q \leq -1$.
\end{proof}
\end{lem}

In next two lemmas, we will show how to generate a homogeneous, $SL(n)$ covariant $L_p$-Minkowski valuation on $\overline{\mathcal{K}}_o ^n$ by a continuous, $SL(n)$ contravariant normalized symmetric $L_p$-Blaschke valuation which is homogeneous of degree $q$ on $\mathcal{K}_o ^n$, where $q \leq -1$.

\begin{lem}\label{Min1}
Let $ Z : \mathcal{K}_o ^n \rightarrow \langle \mathcal{K}_c ^n,\widetilde{\#}_p \rangle$ be a continuous, $SL(n)$ contravariant valuation
which is homogeneous of degree $q = -1$. Define the map $\overline{Z}_1:~\overline{\mathcal{K}}_o ^n \to \langle \overline{\mathcal{K}}_o ^n,+_p \rangle$ by
\begin{align*}
    h(\overline{Z}_1 K,x)^p=\begin{cases}
    C_p \frac{{S_p(ZK,\cdot)}}{V(ZK)}(x) & \dim K =n, \\
    C_p \frac{{S_p(Z[K,\pm b_{k+1},\cdots,\pm b_n],\cdot)}}{V(Z[K,\pm b_{k+1},\cdots,\pm b_n])}(\pi _K x) & \dim K = k<n,
    \end{cases}
\end{align*}
for every $x \in \mathbb{R}^n$, where the $b_{k+1},\cdots,b_n$ are an orthonormal basis of the orthogonal complement of $\text{lin}~K$ and $\pi _K$ is the orthogonal projection onto $\text{lin}~K$.  Then $\overline{Z}_1$ is a $SL(n)$ covariant $L_p$-Minkowski valuation which is homogeneous of degree $1$.
\begin{proof}
In order to show that $\overline{Z}_1$ is well defined, suppose that $\dim K =k<n$ and $b_{k+1},\cdots,b_n$ as well as
$c_{k+1},\cdots,c_n$ are two different orthonormal bases of $(\text{lin}~K)^ \perp$.  Fix an orthonormal basis $b_1,\cdots,b_k$ of $\text{lin}~K$.  Denote
by $\theta$ a proper rotation with $\theta b_i = b_i,~i=1,\cdots,k$ and $\theta b_i \in \{ \pm c_i \},~i=k+1,\cdots,n$.  Then the contravariance of $Z$ and relation (\ref{ncp}) induce that
\begin{align*}
C_p \frac{{S_p(Z[K,\pm c_{k+1},\cdots,\pm c_n],\cdot)}}{V(Z[K,\pm c_{k+1},\cdots,\pm c_n])}(\pi _K x)
&= C_p \frac{{S_p(Z\theta [K,\pm b_{k+1},\cdots,\pm b_n],\cdot)}}{V(Z\theta [K,\pm b_{k+1},\cdots,\pm b_n])}(\pi _K x) \\
&= C_p \frac{{S_p(\theta Z [K,\pm b_{k+1},\cdots,\pm b_n],\cdot)}}{V(\theta Z [K,\pm b_{k+1},\cdots,\pm b_n])}(\pi _K x) \\
&= C_p \frac{{S_p(Z [K,\pm b_{k+1},\cdots,\pm b_n],\cdot)}}{V(Z [K,\pm b_{k+1},\cdots,\pm b_n])}(\theta ^{-1} \pi _K x) \\
&= C_p \frac{{S_p(Z [K,\pm b_{k+1},\cdots,\pm b_n],\cdot)}}{V(Z [K,\pm b_{k+1},\cdots,\pm b_n])}(\pi _K x).
\end{align*}
Thus, $\overline{Z}_1$ is well defined.

Next, we show that $\overline{Z}_1$ is a $L_p$-Minkowski valuation.  Suppose that $K,L \in \overline{\mathcal{K}}_o ^n$ such that $K \cup L \in \overline{\mathcal{K}}_o ^n$ and let $k$ be an integer not larger than n.  If $\dim (K \cup L) = k$, then one of the following four cases is valid:\\
$(1_k)~~~~\dim K = k,\dim L = k, \dim K \cap L =k, ~~~~0 \leq k \leq n,$ \\
$(2_k)~~~~\dim K = k,\dim L = k, \dim K \cap L =k-1, ~~~~~1 \leq k \leq n,$ \\
$(3_k)~~~~\dim K = k,\dim L = k-1, ~~~~~~~~~~~~~~~1 \leq k \leq n,$ \\
$(4_k)~~~~\dim K = k-1,\dim L = k, ~~~~~~~~~~~~~~~1 \leq k \leq n.$ \\
The valuation property trivially holds true for the cases $(3_k)$ and $(4_k)$, since we have $L \subset K$ and
$K \subset L$ respectively in these situations.  Therefore it suffices to prove
\begin{align*}
    h_{\overline{Z}_1(K \cup L)}^p + h_{\overline{Z}_1(K \cap L)}^p = h_{\overline{Z}_1 K}^p +h_{\overline{Z}_1 L}^p
\end{align*}
for the cases $(1_k),0 \leq k \leq n$ and $(2_k),1 \leq k \leq n$.

Let us start with the easy case $(1_n)$.  The valuation property of $Z$ implies
\begin{align*}
    \frac{{S_p(Z(K \cup L),\cdot)}}{V(Z(K \cup L))}+ \frac{{S_p(Z(K \cap L),\cdot)}}{V(Z(K \cap L))}
= \frac{{S_p(ZK,\cdot)}}{V(ZK)} + \frac{{S_p(ZL,\cdot)}}{V(ZL)},
\end{align*}
and thus
\begin{align*}
    C_p \frac{{S_p(Z(K \cup L),\cdot)}}{V(Z(K \cup L))}+ C_p \frac{{S_p(Z(K \cap L),\cdot)}}{V(Z(K \cap L))}
= C_p \frac{{S_p(ZK,\cdot)}}{V(ZK)} + C_p \frac{{S_p(ZL,\cdot)}}{V(ZL)}.
\end{align*}
Hence, the definition of $\overline{Z}_1$ immediately proves the assertion.  Next we deal with the case $(1_k)$, $0 \leq k <n$.
Note that
\begin{align*}
[K,\pm b_{k+1},\cdots,\pm b_n] \cup [L,\pm b_{k+1},\cdots,\pm b_n] = [K \cup L,\pm b_{k+1},\cdots,\pm b_n], \\
[K,\pm b_{k+1},\cdots,\pm b_n] \cap [L,\pm b_{k+1},\cdots,\pm b_n] = [K \cap L,\pm b_{k+1},\cdots,\pm b_n].
\end{align*}
Since $\text{lin}~K = \text{lin}~ L = \text{lin}~ (K \cup L) = \text{lin}~ (K \cap L)$, we have $\pi _K x = \pi _L x = \pi _{(K \cup L)} x
= \pi _{(K \cap L)} x$.  With the valuation property of case $(1_n)$ proved above, we get
\begin{align*}
&C_p \frac{S_p(Z[K \cup L,\pm b_{k+1},\cdots,\pm b_n],\cdot)}{V(Z[K \cup L,\pm b_{k+1},\cdots,\pm b_n])}(x)
+ C_p \frac{S_p(Z[K \cap L,\pm b_{k+1},\cdots,\pm b_n],\cdot)}{V(Z[K \cap L,\pm b_{k+1},\cdots,\pm b_n])}(x) \\
&=C_p \frac{S_p(Z[K,\pm b_{k+1},\cdots,\pm b_n],\cdot)}{V(Z[K,\pm b_{k+1},\cdots,\pm b_n])}(x)
+C_p \frac{S_p(Z[L,\pm b_{k+1},\cdots,\pm b_n],\cdot)}{V(Z[L,\pm b_{k+1},\cdots,\pm b_n])}(x)
\end{align*}
for every $x \in \mathbb{R}^n$. Changing $x$ to $\pi _K x$, we get the positive assertion of the cases $(1_k)$.

Now we consider the case $(2_k),~1 \leq k \leq n$.  It is enough to show
\begin{align}\label{2}
    h_{\overline{Z}_1 K}^p + h_{\overline{Z}_1 (K \cap u^ \perp)}^p
= h_{\overline{Z}_1 (K \cap u^ +)}^p + h_{\overline{Z}_1 (K \cap u^ -)}^p
\end{align}
for $\dim K = k$ and a unit vector $u \in \text{lin}~K$ such that $K \cap u^+,K \cap u^-$ are both $k$-dimensional.  Notice that if $k=n$, then $\pi _K x =x$.  So we will prove the case $(2_k)$ without distinguishing between $k=n$ and $k<n$.
Let $b_1,\cdots,b_n$ be an orthonormal basis of $\mathbb{R}^n$ such that $\text{lin}~K = \text{lin}~ \{ b_1,\cdots,b_k \}$, and $u=b_k$.  With the valuation property of case $(1_k)$ proved above, we have
\begin{align}\label{1}
& C_p \frac{{S_p(Z[K,\pm sb_k,\pm b_{k+1},\cdots,\pm b_n],\cdot)}}{V(Z[K,\pm sb_k,\pm b_{k+1},\cdots,\pm b_n])}(\pi _K x) \nonumber \\
& \qquad +C_p \frac{{S_p(Z[K \cap b_k ^ \perp,\pm sb_k,\pm b_{k+1},\cdots,\pm b_n],\cdot)}}{V(Z[K \cap b_k ^ \perp,\pm sb_k,\pm b_{k+1},\cdots,\pm b_n])}(\pi _K x) \nonumber \\
&=C_p \frac{{S_p(Z[K \cap b_k ^+,\pm sb_k,\pm b_{k+1},\cdots,\pm b_n],\cdot)}}{V(Z[K \cap b_k ^+,\pm sb_k,\pm b_{k+1},\cdots,\pm
b_n])}(\pi _K x) \nonumber \\
& \qquad +C_p \frac{{S_p(Z[K \cap b_k ^-,\pm sb_k,\pm b_{k+1},\cdots,\pm b_n],\cdot)}}{V(Z[K \cap b_k ^-,\pm sb_k,\pm b_{k+1},\cdots,\pm b_n])}(\pi _K x)
\end{align}
for sufficiently small $s>0$.  Define a linear map $\phi$ by
$$\phi b_k = sb_k,
\phi b_i =b_i,~~~~~i=1,\cdots,k-1,k+1,\cdots,n.$$
Note that $\det \phi$ is independent of the choice of orthonormal basis of $\mathbb{R}^n$, so $\det \phi = s$.  The contravariance of $Z$, and relations (\ref{ctho}) as well as (\ref{ncp}) give
\begin{align}\label{2k1}
&C_p \frac{{S_p(Z[K \cap b_k ^ \perp,\pm sb_k,\pm b_{k+1},\cdots,\pm b_n],\cdot)}}{V(Z[K \cap b_k ^ \perp,\pm sb_k,\pm b_{k+1},\cdots,\pm b_n])}(\pi _K x) \nonumber \\
=&C_p \frac{{S_p(Z \phi [K \cap b_k ^ \perp,\pm b_k,\pm b_{k+1},\cdots,\pm b_n],\cdot)}}{V(Z \phi [K \cap b_k ^ \perp,\pm b_k,\pm b_{k+1},\cdots,\pm b_n])}(\pi _K x) \nonumber \\
=&C_p \frac{{S_p(s^\frac{{q+1}}{n} \phi ^{-t} Z[K \cap b_k ^ \perp,\pm b_k,\pm b_{k+1},\cdots,\pm b_n],\cdot)}}{V(s^\frac{{q+1}}{n} \phi ^{-t} Z[K \cap b_k ^ \perp,\pm b_k,\pm b_{k+1},\cdots,\pm b_n])}(\pi _K x) \nonumber \\
=&s^\frac{{-(q+1)p}}{n} C_p \frac{{S_p(Z[K \cap b_k ^ \perp,\pm b_k,\pm b_{k+1},\cdots,\pm b_n],\cdot)}}{V(Z[K \cap b_k ^ \perp,\pm b_k,\pm b_{k+1},\cdots,\pm b_n])}(\phi ^t \pi _K x).
\end{align}
Note that $\mathop {\lim }\limits_{s \to 0^ + } \phi ^t \pi _K x = \pi _{K \cap b_k^ \perp}x$.  Since $q=-1$,
\begin{align*}
&\mathop {\lim }\limits_{s \to 0^ + }C_p \frac{{S_p(Z[K \cap b_k ^ \perp,\pm sb_k,\pm b_{k+1},\cdots,\pm b_n],\cdot)}}{V(Z[K \cap b_k ^ \perp,\pm sb_k,\pm b_{k+1},\cdots,\pm b_n])}(\pi _K x) \nonumber \\
&= C_p \frac{{S_p(Z[K \cap b_k ^ \perp,\pm b_k,\pm b_{k+1},\cdots,\pm b_n],\cdot)}}{V(Z[K \cap b_k ^ \perp,\pm b_k,\pm b_{k+1},\cdots,\pm b_n])}(\pi _{K \cap b_k^ \perp}x).
\end{align*}
So if $s$ tends to zero in (\ref{1}), then we immediately obtain (\ref{2}).
Hence we proved that $\overline{Z}_1$ is a $L_p$-Minkowski valuation.  Moreover, it is easy to calculate that $\overline{Z}_1$ is a $SL(n)$ covariant $L_p$-Minkowski valuation which is homogeneous of degree $1$ on $n$-dimensional convex bodies.
Lemma \ref{Min} implies that $\overline{Z}_1$ is a $SL(n)$ covariant $L_p$-Minkowski valuation which is homogeneous of degree $1$ on $\overline{\mathcal{K}}_o ^n$.
\end{proof}
\end{lem}

\begin{lem}\label{Min2}
Let $ Z : \mathcal{K}_o ^n \rightarrow \langle \mathcal{K}_c ^n,\widetilde{\#}_p \rangle$ be a continuous, $SL(n)$ contravariant valuation
which is homogeneous of degree $q<-1$. Define the map $\overline{Z}_2:~\overline{\mathcal{K}}_o ^n \to \langle \overline{\mathcal{K}}_o ^n,+_p \rangle$ by
\begin{align*}
    h(\overline{Z}_2 K,x)^p=\begin{cases}
    C_p \frac{{S_p(ZK,\cdot)}}{V(ZK)}(x) & \dim K =n, \\
    0 & \dim K = k<n,
    \end{cases}
\end{align*}
for every $x \in \mathbb{R}^n$.  Then $\overline{Z}_2$ is a $SL(n)$ covariant $L_p$-Minkowski valuation which is homogeneous of degree $r = -q$.
\begin{proof}
We use the notation of Lemma \ref{Min1}.  Since the case $(1_n)$ is the same as in Lemma \ref{Min1}, and the cases $(1_k),~0 \leq k < n$, $(2_k),~1 \leq k < n$ are trivially true, we just need to consider the case $(2_n)$.

Hence we need to show
\begin{align}\label{3}
    h_{\overline{Z}_2 K}^p + h_{\overline{Z}_2 (K \cap u^ \perp)}^p
= h_{\overline{Z}_2 (K \cap u^ +)}^p + h_{\overline{Z}_2 (K \cap u^ -)}^p
\end{align}
for $\dim K = n$ and a unit vector $u \in \mathbb{R}^n$ such that $K \cap u^+,K \cap u^-$ are both $n$-dimensional.  Let $b_1,\cdots,b_n$ be an orthonormal basis of $\mathbb{R}^n$ such that $u=b_n$.
Comparing with the proof of Lemma \ref{Min1}, we just need to show the relation (\ref{2k1}) of the case $k = n$ tends to zero for $q<-1$ when $s$ tends to zero.  Actually, the relation (\ref{2k1}) of the case $k = n$ is
\begin{align*}
C_p \frac{{S_p(Z[K \cap b_n ^ \perp,\pm sb_n],\cdot)}}{V(Z[K \cap b_n ^ \perp,\pm sb_n])}(x) =
s^{\frac{{-(q+1)p}}{n}}C_p \frac{S_p(Z[K \cap b_n ^ \perp,\pm b_n],\cdot)}{V(Z[K \cap b_n ^ \perp,\pm b_n])}(\phi ^t x),
\end{align*}
where $\phi$ is a linear map defined by $\phi b_n = sb_n,\phi b_i =b_i,~~~~~i=1,\cdots,n-1$.
Since $q<-1$,
\begin{align*}
\mathop {\lim }\limits_{s \to 0^ + }C_p \frac{{S_p(Z[K \cap b_n ^ \perp,\pm sb_n], \cdot)}}{V(Z[K \cap b_n ^ \perp,\pm sb_n])}(x)=0.
\end{align*}
Hence, $\overline{Z}_2$ is a $L_p$-Minkowski valuation. Moreover, it is easy to calculate that $\overline{Z}_2$ is a $SL(n)$ covariant $L_p$-Minkowski valuation which is homogeneous of degree $r = -q $.
\end{proof}
\end{lem}

For $p>1$, the following lemma shows that every support set of a $L_p$-projection body consists of precisely one point.  It will help to rule out the existence of continuous, normalized symmetric $L_p$-Blaschke valuations which are homogeneous of degree $-1$ (see Theorem \ref{npbt} and Theorem \ref{ncop} for more details).  A similar result for $p=1$ can be found in Schneider \cite[Lemma 3.5.5]{Schn2}.

For $K \in \mathcal{K}^n,~e \in S^{n-1}$, write $K_e := \{x \in K|x \cdot e = h(K,e) \}$.
\begin{lem}\label{pgz}
For $p>1$, if the support function of the convex body $K \in \mathcal{K}^n$ is given by
$$h(K,u)= (\int _{S^{n-1}} |u \cdot v|^p d \mu (v))^{1/p}$$
for $u \in S^{n-1}$, with an even signed measure $\mu$, then, for $e \in S^{n-1}$,
$$h(K_e,u)=v_e \cdot u$$
for $u \in S^{n-1}$, where $v_e = 2 (\int _{S^{n-1}} |e \cdot v|^p d \mu (v))^{\frac{1}{p}-1} \int _{e^+} (e \cdot v)^{p-1} v d \mu (v)$.
\begin{proof}
The assertion of the lemma is true for $u= \pm e$, since $h(K_e,\pm e)=\pm h(K,e)$.  Hence we may assume that $u$ and $e$ are linearly independent. Note that $h(K_e,u)= \mathop {\lim }\limits_{s \to 0^ + } \frac{h(K,e+su)-h(K,e)}{s}$ (see Schneider \cite[Theorem 1.7.2]{Schn2}).  Put
\begin{align*}
&A_s:=\{v \in S^{n-1}|e \cdot v >0, (e + su) \cdot v >0 \}, \\
&B_s:=\{v \in S^{n-1}|e \cdot v \leq 0, (e + su) \cdot v >0 \}, \\
&C_s:=\{v \in S^{n-1}|e \cdot v >0, (e + su) \cdot v \leq 0 \}.
\end{align*}
We obtain
\begin{align*}
&h(K_e,u)= \mathop {\lim }\limits_{s \to 0^ + } \frac{h(K,e+su)-h(K,e)}{s} \\
&= \frac{1}{p}(\int _{S^{n-1}} |e \cdot v|^p d \mu (v))^{\frac{1}{p}-1} \mathop {\lim }\limits_{s \to 0^ + } \frac{1}{s}(\int _{S^{n-1}} |(e+su) \cdot v|^p d \mu (v)-\int _{S^{n-1}} |e \cdot v|^p d \mu (v)),
\end{align*}
and
\begin{align*}
&\mathop {\lim }\limits_{s \to 0^ + } \frac{1}{s}(\int _{S^{n-1}} |(e+su) \cdot v|^p d \mu (v)-\int _{S^{n-1}} |e \cdot v|^p d \mu (v)) \\
&=2 \mathop {\lim }\limits_{s \to 0^ + } \frac{1}{s}(\int _{A_s \cup B_s} ((e+su) \cdot v)^p d \mu (v)-\int _{A_s \cup C_s} (e \cdot v)^p d \mu (v)) \\
&=2p \mathop {\lim }\limits_{s \to 0^ + } \int _{A_s \cup B_s} (e \cdot v)^{p-1} (u \cdot v) d \mu (v) \\
& \qquad + \mathop {\lim }\limits_{s \to 0^ + } \int _{A_s \cup B_s} p(p-1)(e \cdot v)^{p-2} (u \cdot v)^2s + o(s)d \mu (v) \\
& \qquad +2 \mathop {\lim }\limits_{s \to 0^ + } \frac{1}{s}\int _{B_s} (e \cdot v)^p d \mu (v)
- 2 \mathop {\lim }\limits_{s \to 0^ + } \frac{1}{s} \int _{C_s} (e \cdot v)^p d \mu (v).
\end{align*}
Let
$$\mu _+ (E) = \sup \{ \mu (A)|A \subset E~\rm{and}~A~\rm{is~a~Borel~set~of}~S^{n-1}\},$$ $$\mu _- (E) = - \inf \{ \mu (A)|A \subset E~\rm{and}~A~\rm{is~a~Borel~set~of}~S^{n-1}\},$$ $$\mu '(E) = \mu _+(E) + \mu_-(E)$$
for every Borel set $E$ of $S^{n-1}$. We get
\begin{align*}
&\big| \int _{A_s \cup B_s} p(p-1)(e \cdot v)^{p-2} (u \cdot v)^2s + o(s)d \mu (v) \big| \\
& \qquad \qquad \leq \int _{S^{n-1}} \big| p(p-1)(e \cdot v)^{p-2} (u \cdot v)^2s + o(s) \big| d \mu ' (v) \xrightarrow{{s \to 0^ +  }} 0.
\end{align*}
For $v \in B_s$, we have $|e \cdot v| \leq cs$ with a constant $c$ independent of $s$. Writing
$$B_s':=\{v \in S^{n-1} | e \cdot v < 0, (e + su) \cdot v >0 \},$$
we obtain
\begin{align*}
\big| \frac{1}{s}\int _{B_s} (e \cdot v)^p d \mu (v) \big|= \big| \frac{1}{s} \int _{B_s'} (e \cdot v)^p d \mu (v) \big| \leq c^p s^{p-1} \mu (B_s').
\end{align*}
Since (in the set-theoretic sense) $\mathop {\lim }\limits_{s \to 0^ + } B_s'=\emptyset$, we have $\mathop{\lim }\limits_{s \to 0^ + } \mu ' (B_s')=0$. With $p>1$, we get
\begin{align*}
\mathop {\lim }\limits_{s \to 0^ + } \frac{1}{s}\int _{B_s} (e \cdot v)^p d \mu (v)=0.
\end{align*}
From $\mathop {\lim }\limits_{s \to 0^ + } C_s=\emptyset$, we similarly find
\begin{align*}
\mathop {\lim }\limits_{s \to 0^ + } \frac{1}{s}\int _{C_s} (e \cdot v)^p d \mu (v)=0.
\end{align*}
Further, $\mathop {\lim }\limits_{s \to 0^ + } A_s = e^+ \backslash e^ \perp$, $\mathop {\lim }\limits_{s \to 0^ + } B_s = \{v \in S^{n-1} | e \cdot v = 0, u \cdot v >0\}$, and $p>1$, we get
\begin{align*}
\mathop {\lim }\limits_{s \to 0^ + } \int _{A_s \cup B_s} (e \cdot v)^{p-1} (u \cdot v) d \mu (v) = \int _{e^+} (e \cdot v)^{p-1} (u \cdot v) d \mu (v).
\end{align*}
Finally, we get
\begin{align*}
h(K_e,u) &= 2 (\int _{S^{n-1}} |e \cdot v|^p d \mu (v))^{\frac{1}{p}-1} \int _{e^+} (e \cdot v)^{p-1} (u \cdot v) d \mu (v) \\
&= \big(2 (\int _{S^{n-1}} |e \cdot v|^p d \mu (v))^{\frac{1}{p}-1} \int _{e^+} (e \cdot v)^{p-1} v d \mu (v) \big) \cdot u,
\end{align*}
which completes the proof of the lemma.
\end{proof}
\end{lem}

To classify continuous, homogeneous, $SL(n)$ contravariant normalized symmetric $L_p$-Blaschke valuations, we need the following results from Ludwig \cite{ludwigMin}.

For $-1 \leq \tau \leq 1$, define $M _p ^ \tau :\overline{\mathcal{K}}_o ^n \to \overline{\mathcal{K}}_o ^n$ by
\begin{align*}
    h^p(M _p ^ \tau K,v)=\int_K (|v \cdot x| + \tau (v \cdot x))^p dx
\end{align*}
for $v \in \mathbb{R}^n$.  In particular, $M _p ^ 0 K$ is a dilate of the $L_p$-centroid body, if $V(K)>0$.

A polytope is the convex hull of finitely many points in $\mathbb{R}^n$. Let $\mathcal{P}_o ^n$ be the set of $n$-dimensional polytopes which contain the origin, $\overline{\mathcal{P}}_o ^n$ the set of polytopes which contain the origin. Let $\xi _o(P)$ denote the set of edges of a polytope $P$ which contain the origin.

\begin{lem}\label{ludwigtrp}
\rm{\cite{ludwigMin}} Let $ Z : \overline{\mathcal{P}}_o ^n \rightarrow \langle \overline{\mathcal{K}}_o ^n,+_p \rangle$, $n \geq 3$, be a $L_p$-Minkowski valuation, $p > 1$, which is $SL(n)$
covariant and homogeneous of degree $r$.  If $r = n/p +1$, then there are constants $ a \geq 0 $ and $-1 \leq \tau \leq 1 $ such
that $$Z P =a M _p ^ \tau P $$
for every $ P \in \overline{\mathcal{P}}_o ^n $.  If $r = 1$, then there are constants $ a,b \geq 0 $ such
that $$Z P =aP +_p b(-P) $$
for every $ P \in \overline{\mathcal{P}}_o ^n $.  In all other cases, $Z P = \{o\} $ for every $ P \in \overline{\mathcal{P}}_o ^n $.

Let $ Z : \overline{\mathcal{P}}_o ^2 \rightarrow \langle \overline{\mathcal{K}}_o ^2,+_p \rangle$, be a $L_p$-Minkowski valuation, $p > 1$, which is $SL(2)$ covariant and homogeneous of degree $r$.  If $r = 2/p +1$, then there are constants $ a \geq 0 $ and $-1 \leq \tau \leq 1 $ such
that $$Z P =a M _p ^ \tau P $$
for every $ P \in \overline{\mathcal{P}}_o ^2 $.  If $r = 1$, then there are constants $ a_0,b_0 \geq 0 $ and $a_i,b_i \in \mathbb{R}$ with $a_0+a_i,~b_0+b_i \geq 0,~i=1,2$ such that
\begin{align*}
Z P =a_0P +_p b_0(-P)+_p \sum \nolimits ^p (a_i E_i +_p b_i (-E_i))
\end{align*}
for every $ P \in \overline{\mathcal{P}}_o ^2 $, where $\sum \nolimits ^p$ denotes the $L_p$-Minkowski sum, and the sum is taken over $E_i \in \xi _o (P)$.  If $r=2/p-1$, then there are constants $a \geq 0$ and $-1 \leq \tau \leq 1$ such that
$$ZP= a \psi _{\pi /2} \hat{\Pi} _p ^ \tau P$$
for every $ P \in \overline{\mathcal{P}}_o ^2 $. Here $\hat{\Pi} _p ^ \tau P$ is defined by the relation (\ref{pp}).  In all other cases, $Z P = \{o\} $for every $ P \in \overline{\mathcal{P}}_o ^2 $.
\end{lem}

Now we can classify continuous, homogeneous, $SL(n)$ contravariant normalized symmetric $L_p$-Blaschke valuations.
\begin{thm}\label{npbt}
    Let $n \geq 2$, $p > 1$ and $p$ not an even integer. If $Z : \mathcal{K}_o ^n \rightarrow \langle \mathcal{K}_c ^n,\widetilde{\#}_p \rangle$ is a continuous, homogeneous, $SL(n)$ contravariant valuation, then there exists a constant $c>0$ such that
$$ZK=c \widetilde{\Lambda} _c ^p K$$
for every $K \in \mathcal{K}_o ^n$.
\begin{proof}
Let $q$ be the degree of homogeneity of $Z$. Lemma \ref{3.1} shows that $q \leq -1$.

If $q=-1$, then $\overline{Z}_1$, introduced in Lemma \ref{Min1}, is a $SL(n)$ covariant $L_p$-Minkowski valuation which is homogeneous of degree $1$.  If $n \geq 3$, from Lemma \ref{ludwigtrp}, we derive that there are constants $ a,b \geq 0 $ such that
\begin{align*}
\overline{Z}_1 P =aP +_p b(-P)
\end{align*}
for every $ P \in \overline{\mathcal{P}}_o ^n $.  If $n=2$, from Lemma \ref{ludwigtrp}, we derive that there are constants $ a_0,b_0 \geq 0 $ and $a_i,b_i \in \mathbb{R}$ with $a_0+a_i,~b_0+b_i \geq 0,~i=1,2$ such that
\begin{align*}
\overline{Z}_1 P =a_0P +_p b_0(-P)+_p \sum \nolimits ^p (a_i E_i +_p b_i (-E_i))
\end{align*}
for every $ P \in \overline{\mathcal{P}}_o ^n$, where the sum is taken over $E_i \in \xi _o (P)$.  For $P_0 = [\pm e_1, \cdots, \pm e_n]$, we have
\begin{align*}
\frac{{\Pi _p ZP_0}}{V(ZP_0)^{1/p}} = cP_0,
\end{align*}
for a suitable $c \geq 0$ when $n \geq 2$.  Assumption that $Z$ does not contain $\{o \}$ in its range gives $c>0$.  For $p>1$, every support set of a $L_p$-projection body consists of precisely one point (Lemma \ref{pgz}).  However, $P_0$ has the support set $[e_1,\cdots,e_n]$ which does not consist of precisely one point.  That is a contradiction.

If $q = -n/p-1$, then $\overline{Z}_2$, introduced in Lemma \ref{Min2}, is a $SL(n)$ covariant $L_p$-Minkowski valuation which is homogeneous of degree $n/p + 1$.  For $n \geq 2$, from Lemma \ref{ludwigtrp}, we infer that the existence of constants $ a \geq 0 $ and $-1 \leq \tau \leq 1 $ such that
\begin{align*}
\overline{Z}_2 P =a M _p ^ \tau P
\end{align*}
for every $ P \in \overline{\mathcal{P}}_o ^n $.  Assumption that $Z$ does not contain $\{o \}$ in its range gives $a>0$.  Since $\overline{Z}_2 P$ is origin-symmetric, we deduce that $\tau = 0$.  Thus,
$\frac{\Pi _p ZP}{V(ZP)^{1/p}} =a M_p ^0 P$ for every $P \in \mathcal{P}_o ^n$.  Since the operators $\frac{\Pi _p Z}{V^{1/p}}$ and $\Gamma _p$ are continuous on $\mathcal{K}_o ^n$, and $\mathcal{P}_o ^n$ is dense in $\mathcal{K}_o ^n$, we obtain
\begin{align*}
\frac{{\Pi _p ZK}}{V(ZK)^{1/p}} =a M_p ^0 K
\end{align*}
for every $K \in \mathcal{K}_o ^n$.  By rewriting this in terms of the $L_p$-cosine transforms (via relation (\ref{2.11}) and $(c_{n,p}V(K))^\frac{1}{p}\Gamma _p K = M_p ^0 K$), we get
\begin{align*}
C_p \frac{{S_p(ZK,\cdot)}}{V(ZK)} = b C_p (\rho _K (\cdot)^{n+p}) = b C_p (\frac{{1}}{2} \rho _K (\cdot)^{n+p} + \frac{{1}}{2} \rho _{-K} (\cdot)^{n+p})
\end{align*}
for a suitable constant $b > 0$.  Since $S_p(ZK,\cdot)$ is an even measure, the injectivity property (\ref{inj}) and the definition of the normalized symmetric $L_p$-curvature image operator finally shows
\begin{align}
ZK = c \widetilde{\Lambda} _c ^p K
\end{align}
for a suitable constant $c>0$.

If $q=-2/p +1$ and $n=2$, then $\overline{Z}_2$, introduced in Lemma \ref{Min2}, is a $SL(n)$ covariant $L_p$-Minkowski valuation which is homogeneous of degree $2/p -1$.  By Lemma \ref{ludwigtrp}, there are constants $ a \geq 0 $ and $-1 \leq \tau \leq 1 $ such
that $$\overline{Z}_2 P =a \hat{\Pi} _p ^ \tau P $$
for every $ P \in \overline{\mathcal{P}}_o ^n $.  $\hat{\Pi} _p ^ \tau$ is not continuous on $\mathcal{P}_o ^n$ while $\frac{\Pi _p Z}{V^{1/p}}$ is continuous on $\mathcal{P}_o ^n$. Thus, that is a contradiction.

In all other cases, $\overline{Z}_2$, introduced in Lemma \ref{Min2}, is a $SL(n)$ covariant $L_p$-Minkowski valuation which is homogeneous of degree $r$, where $r \neq 1$, $r \neq n/p+1$ for $n \geq 2$ and $r \neq 2/p - 1$ as an addition for $n=2$.  By Lemma \ref{ludwigtrp}, $\overline{Z}_2 P=\{o\}$ for every $P \in \overline{\mathcal{P}}_o ^n$.  So
\begin{align}
h_{\overline{Z}_2 P}(x)^p = C_p \frac{{S_p(ZP,\cdot)}}{V(ZP)}(x) = 0
\end{align}
for every $P \in \mathcal{P}_o ^n$.  $S_p(ZP,\cdot)$ is an even measure since $ZP$ is an origin-symmetric convex body.  Thus, by relation (\ref{inj}), we have $S_p(ZP,\cdot)=0$.  That is a contradiction.
\end{proof}
\end{thm}

Hence, Theorem \ref{th3.1} and Theorem \ref{npbt} directly imply Theorem \ref{1.1}.
\subsection{ The covariant case }
The following Lemma \ref{4.1}, Lemma \ref{coMin1} and Lemma \ref{coMin2} are the counterparts of Lemma \ref{3.1}, Lemma \ref{Min1} and Lemma \ref{Min2} respectively.
\begin{lem}\label{4.1}
If $ Z : \mathcal{K}_o ^n \rightarrow \langle \mathcal{K}_c ^n,\widetilde{\#}_p \rangle$ is a continuous, $SL(n)$ covariant valuation which is homogeneous of degree $q$, then $q \leq -n+1$.
\begin{proof}
Suppose $K \in \mathcal{K}_o ^n$ and $s>0$.  As in the proof of Lemma \ref{3.1}, we get that
\begin{align}\label{co1}
    & \quad C_p \frac{{S_p(Z[K \cap e_n ^ \perp ,\pm se_n],\cdot)}}{V(Z[K \cap e_n ^ \perp ,\pm se_n])}(e_n) \hfill \nonumber \\
    &= C_p \frac{{S_p(Z[K \cap e_n ^+ ,\pm se_n],\cdot)}}{V(Z[K \cap e_n ^+ ,\pm se_n])}(e_n) + C_p \frac{{S_p(Z[K \cap e_n ^- ,\pm se_n],\cdot)}}{V(Z[K \cap e_n ^- ,\pm se_n])}(e_n) \nonumber \\
    & \quad - C_p \frac{{S_p(Z[K,\pm se_n],\cdot)}}{V(Z[K,\pm se_n])}(e_n),
\end{align}
and thus $C_p \frac{{S_p(Z[K \cap e_n ^ \perp ,\pm se_n],\cdot)}}{V(Z[K \cap e_n ^ \perp ,\pm se_n])}(e_n)$ must converge to a finite number as $s \to 0^+$.  (The difference between relation (\ref{ct1}) and relation (\ref{co1}) is that the independent variable of the $L_p$-cosine transform is changed from $e_1$ to $e_n$.)
Define the linear map $\phi$ as before by
$$\phi e_i = e_i, i=1,\cdots,n-1, \phi e_n =se_n.$$
From the $SL(n)$ covariance and homogeneity of $Z$ as well as relation (\ref{coho}) and (\ref{ncp}), we get
\begin{align*}
    C_p \frac{{S_p(Z[K \cap e_n ^ \perp ,\pm se_n],\cdot)}}{V(Z[K \cap e_n ^ \perp ,\pm se_n])}(e_n)
& = C_p \frac{{S_p(s^{ \frac{{q-1}}{n} } \phi Z[K \cap e_n ^ \perp ,\pm e_n],\cdot)}}{V(s^{ \frac{{q-1}}{n} } \phi Z[K \cap e_n ^ \perp ,\pm e_n])}(e_n) \hfill
    \\ & =s^ \frac{{-(q-1)p}}{n} C_p \frac{{S_p(Z[K \cap e_n ^ \perp ,\pm e_n],\cdot)}}{V(Z[K \cap e_n ^ \perp ,\pm e_n])}(\phi ^{-1} e_n). \hfill
\end{align*}
Since $|e_n \cdot u|>0$ for all $u \in S^{n-1} \setminus {e_n ^ \perp}$ and the $L_p$-surface area measure of $n$-dimensional bodies is not concentrated on any great sphere,
we conclude that
\begin{align*}
   & s^ p C_p \frac{{S_p(Z[K \cap e_n ^ \perp ,\pm e_n],\cdot)}}{V(Z[K \cap e_n ^ \perp ,\pm e_n])}(\phi ^{-1} e_n) \\
&= \frac{{1}}{V(Z[K \cap e_n ^ \perp ,\pm e_n])} \int _{S^{n-1}} |e_n \cdot u|^p dS_p(Z[K \cap e_n ^ \perp ,\pm e_n],u)>0.
\end{align*}
Thus, $\frac{{-(q-1)p}}{n} - p \geq 0$, so $q \leq -n+1$.
\end{proof}
\end{lem}

\begin{lem}\label{coMin1}
Let $ Z : \mathcal{K}_o ^n \rightarrow \langle \mathcal{K}_c ^n,\widetilde{\#}_p \rangle$ be a continuous, $SL(n)$ covariant valuation
which is homogeneous of degree $q = -n+1$. Define the map $\overline{Z}_1:~\overline{\mathcal{K}}_o ^n \to \langle \mathcal{K}^n,+_p \rangle$ by
\begin{align*}
    h(\overline{Z}_1 K,x)^p=\begin{cases}
    C_p \frac{{S_p(ZK,\cdot)}}{V(ZK)}(x) & \dim K =n, \\
    C_p \frac{{S_p(Z[K,\pm v],\cdot)}}{V(Z[K,\pm v])}((x \cdot v)v) & \dim K = n-1, \\
    0 & \dim K \leq n-2,
    \end{cases}
\end{align*}
for every $x \in \mathbb{R}^n$, where $v$ is a unit vector perpendicular to $\text{lin}~K$.  Then $\overline{Z}_1$ is a $SL(n)$ contravariant $L_p$-Minkowski valuation which is homogeneous of degree $n-1$.
\begin{proof}
Obviously, the definition of $\overline{Z}_1$ is independent of the choice of $v$, so it is well defined. Next, we show that $\overline{Z}_1$ is a $L_p$-Minkowski valuation.  We still use the notation of the proof of Lemma \ref{Min1}.  The case $(1_n)$ is the same as and the case $(1_{n-1})$ is similar to (Change $\pi _K x$ to $(x \cdot v)v$) the corresponding parts in the proof of Lemma \ref{Min1}.  The cases $(1_k),0 \leq k \leq n-2$ and $(2_k),1 \leq k \leq n-2$ are trivial.

Now we consider the case $(2_n)$.  It is enough to show
\begin{align}\label{co2}
    h_{\overline{Z}_1 K}^p + h_{\overline{Z}_1 (K \cap u^ \perp)}^p
= h_{\overline{Z}_1 (K \cap u^ +)}^p + h_{\overline{Z}_1 (K \cap u^ -)}^p
\end{align}
for $\dim K = n$ and a unit vector $u \in \text{lin}~K$ such that $K \cap u^+,K \cap u^-$ are both $n$-dimensional.  Let $b_1,\cdots,b_n$ be an orthonormal basis of $\mathbb{R}^n$ such that $u=b_n$.  With the valuation property of case $(1_n)$, we have
\begin{align}\label{co11}
&C_p \frac{{S_p(Z[K,\pm sb_n],\cdot)}}{V(Z[K,\pm sb_n])}(x) + C_p \frac{{S_p(Z[K \cap b_n ^ \perp,\pm sb_n],\cdot)}}{V(Z[K \cap b_n ^ \perp,\pm sb_n])}(x) \nonumber \\
&=C_p \frac{{S_p(Z[K \cap b_n ^+,\pm sb_n],\cdot)}}{V(Z[K \cap b_n ^+,\pm sb_n])}(x) +C_p \frac{{S_p(Z[K \cap b_n ^-,\pm sb_n],\cdot)}}{V(Z[K \cap b_n ^-,\pm sb_n])}(x)
\end{align}
for sufficiently small $s>0$.  Define a linear map $\phi$ by
$$\phi b_n = sb_n,
\phi b_i =b_i,~~~~~i=1,\cdots,n-1.$$
The covariance of $Z$, and relations(\ref{coho}) as well as (\ref{ncp}) give
\begin{align}\label{co5}
C_p \frac{{S_p(Z[K \cap b_n ^ \perp,\pm sb_n],\cdot)}}{V(Z[K \cap b_n ^ \perp,\pm sb_n])}(x)
&= s^\frac{{-(q-1)p}}{n} C_p \frac{{S_p(Z[K \cap b_n ^ \perp,\pm b_n],\cdot)}}{V(Z[K \cap b_n ^ \perp,\pm b_n])}(\phi ^{-1} x) \nonumber \\
&= s^{\frac{{-(q-1)p}}{n}-p} C_p \frac{{S_p(Z[K \cap b_n ^ \perp,\pm b_n],\cdot)}}{V(Z[K \cap b_n ^ \perp,\pm b_n])}(s \phi ^{-1} x).
\end{align}
Note that $\mathop {\lim }\limits_{s \to 0^ + } s \phi ^{-1} x = (x \cdot b_n)b_n$.  Since $q=-n+1$,
\begin{align*}
\mathop {\lim }\limits_{s \to 0^ + }C_p \frac{{S_p(Z[K \cap b_n ^ \perp,\pm sb_n],\cdot)}}{V(Z[K \cap b_n ^ \perp,\pm sb_n])}(x)
= C_p \frac{{S_p(Z[K \cap b_n ^ \perp,\pm b_n],\cdot)}}{V(Z[K \cap b_n ^ \perp,\pm b_n])}((x \cdot b_n)b_n).
\end{align*}
So if $s$ tends to zero in (\ref{co11}), then we immediately obtain (\ref{co2}).

The case $(2_{n-1})$ is similar to the case $(2_n)$.  We will show the relation (\ref{co2}) is still true for $\dim K = n-1$ and a unit vector $u \in \text{lin}~K$ such that $K \cap u^+,K \cap u^-$ are both $(n-1)$-dimensional.  Let $b_1,\cdots,b_n$ be an orthonormal basis of $\mathbb{R}^n$ such that $\text{lin}~K = \text{lin}~ \{ b_1,\cdots,b_{n-1} \}$, and $u=b_{n-1}$. Thus choose $v= b_n$.  With the valuation property of case $(1_{n-1})$, we have
\begin{align}\label{co111}
&C_p \frac{{S_p(Z[K,\pm sb_{n-1},\pm b_n],\cdot)}}{V(Z[K,\pm sb_{n-1},\pm b_n])}((x \cdot b_n)b_n) \nonumber \\
& \qquad + C_p \frac{{S_p(Z[K \cap b_{n-1} ^ \perp,\pm sb_{n-1},\pm b_n],\cdot)}}{V(Z[K \cap b_{n-1} ^ \perp,\pm sb_{n-1},\pm b_n])}((x \cdot b_n)b_n) \nonumber \\
=&C_p \frac{{S_p(Z[K \cap b_{n-1} ^+,\pm sb_{n-1},\pm b_n],\cdot)}}{V(Z[K \cap b_{n-1} ^+,\pm sb_{n-1},\pm b_n])}((x \cdot b_n)b_n) \nonumber \\
& \qquad + C_p \frac{{S_p(Z[K \cap b_{n-1} ^-,\pm sb_{n-1},\pm b_n],\cdot)}}{V(Z[K \cap b_{n-1} ^-,\pm sb_{n-1},\pm b_n])}((x \cdot b_n)b_n)
\end{align}
for sufficiently small $s>0$.  Define a linear map $\phi$ by
$$\phi b_{n-1} = sb_{n-1},
\phi b_i =b_i,~~~~~i \neq n-1.$$
The covariance of $Z$, and relations(\ref{coho}) as well as (\ref{ncp}) give
\begin{align*}
&\frac{{S_p(Z[K \cap b_{n-1} ^ \perp,\pm sb_{n-1},\pm b_n],\cdot)}}{V(Z[K \cap b_{n-1} ^ \perp,\pm sb_{n-1},\pm b_n])}((x \cdot b_n)b_n) \\
&= s^\frac{{-(q-1)p}}{n} C_p \frac{{S_p(Z[K \cap b_n ^ \perp,\pm b_n],\cdot)}}{V(Z[K \cap b_n ^ \perp,\pm b_n])}(\phi ^{-1} (x \cdot b_n)b_n).
\end{align*}
Note that $\mathop {\lim }\limits_{s \to 0^ + } \phi ^{-1} (x \cdot b_n)b_n = (x \cdot b_n)b_n$.  Since $q=-n+1$,
\begin{align*}
&\mathop {\lim }\limits_{s \to 0^ + }C_p \frac{{S_p(Z[K \cap b_{n-1} ^ \perp,\pm sb_{n-1},\pm b_n],\cdot)}}{V(Z[K \cap b_{n-1} ^ \perp,\pm sb_{n-1},\pm b_n])}((x \cdot b_n)b_n) =0. \nonumber \\
\end{align*}
So if $s$ tends to zero in (\ref{co111}), then we immediately obtain (\ref{co2}).
Hence we proved that $\overline{Z}_1$ is a $L_p$-Minkowski valuation.

Moreover it is easy to calculate that $\overline{Z}_1$ is a $SL(n)$ contravariant $L_p$-Minkowski valuation which is homogeneous of degree $n-1$ on $n$-dimensional convex bodies.
Lemma \ref{Min} implies that $\overline{Z}_1$ is a $SL(n)$ contravariant $L_p$-Minkowski valuation which is homogeneous of degree $n-1$.
\end{proof}
\end{lem}

\begin{lem}\label{coMin2}
Let $ Z : \mathcal{K}_o ^n \rightarrow \langle \mathcal{K}_c ^n,\widetilde{\#}_p \rangle$ be a continuous, $SL(n)$ covariant valuation
which is homogeneous of degree $q<-n+1$. Define the map $\overline{Z}_2:~\overline{\mathcal{K}}_o ^n \to \langle \mathcal{K}^n,+_p \rangle$ by
\begin{align*}
    h(\overline{Z}_2 K,x)^p=\begin{cases}
    C_p \frac{{S_p(ZK,\cdot)}}{V(ZK)}(x) & \dim K =n, \\
    0 & \dim K = k<n,
    \end{cases}
\end{align*}
for every $x \in \mathbb{R}^n$.  Then $\overline{Z}_2$ is a $SL(n)$ contravariant $L_p$-Minkowski valuation which is homogeneous of degree $r = -q$.
\begin{proof}
To prove that $\overline{Z}_2$ is a $L_p$-Minkowski valuation, as in the proof of Lemma \ref{Min2}, we just need to show
\begin{align}
\mathop {\lim }\limits_{s \to 0^ + }C_p \frac{{S_p(Z[K \cap b_n ^ \perp,\pm sb_n],\cdot)}}{V(Z[K \cap b_n ^ \perp,\pm sb_n])}(x)
= 0.
\end{align}
Actually, since $q<-n+1$, by the relation (\ref{co5}), we immediately get the conclusion.

Moreover, it is easy to calculate that $\overline{Z}_2$ is a $SL(n)$ covariant $L_p$-Minkowski valuation which is homogeneous of degree $r = -q$.
\end{proof}
\end{lem}

As in the contravariant case, we also need following results from \cite{ludwigMin} to classify $SL(n)$ covariant normalized symmetric $L_p$-Blaschke valuations.

For $-1 \leq \tau \leq 1$, define $\Pi _p ^ \tau$ on the set of all convex bodies containing the origin in their interiors by
\begin{align*}
    h(\Pi _p ^ \tau K,v)^p= \int_{S^{n-1}} (|v \cdot u| + \tau (v \cdot u))^p dS_p(K,u)
\end{align*}
for $v \in \mathbb{R}^n$.  In particular, $\Pi _p ^ 0 K$ is the $L_p$-projection body of $K$.  To extend the operator $\Pi _p ^ \tau$ to polytopes that contain the origin in their boundaries, for $P \in \overline{\mathcal{P}}_o ^n$, the set of polytopes which contain the origin, define $\hat{\Pi} _p ^ \tau P$ by
\begin{align}\label{pp}
    h(\hat{\Pi} _p ^ \tau P,v)^p= \int _{S^{n-1} \setminus \omega _o (P)} (|v \cdot u| + \tau (v \cdot u))^p dS_p(P,u),
\end{align}
where $\omega _o (P)$ is the set of outer unit normal vectors to facets of $P$ that contain the origin.

\begin{lem}\label{ludwigcop} \rm{\cite{ludwigMin}}
Let $ Z : \overline{\mathcal{P}}_o ^n \rightarrow \langle \mathcal{K}^n,+_p \rangle$ be a $L_p$-Minkowski valuation, $p > 1,n \geq 3$, which is $SL(n)$ contravariant and homogeneous of degree $r$.  If $r = n/p-1$, then there are constants $ a \geq 0 $ and $-1 \leq \tau \leq 1 $ such that $$Z P =a \hat{\Pi} _p ^ \tau P $$
for every $ P \in \overline{\mathcal{P}}_o ^n $.  In all other cases, $Z P = \{o\} $ for every $ P \in \overline{\mathcal{P}}_o ^n $.

Let $ Z : \overline{\mathcal{P}}_o ^2 \rightarrow \langle \mathcal{K}^2,+_p \rangle$ be a $L_p$-Minkowski valuation, $p > 1$, which is $SL(2)$ contravariant and homogeneous of degree $r$.  If $r=2/p+1$, then there are constants $ a \geq 0 $ and $-1 \leq \tau \leq 1 $ such that
$$Z P =a \psi _{\pi /2} M _p ^ \tau P $$
for every $ P \in \overline{\mathcal{P}}_o ^2 $.  If $r=1$, then there are constants $ a_0,b_0 \geq 0 $ and $a_i,b_i \in \mathbb{R}$ with $a_0+a_i,~b_0+b_i \geq 0,~i=1,2$ such that
\begin{align*}
Z P =\psi _{\pi /2} (a_0P +_p b_0(-P)+_p \sum \nolimits ^p (a_i E_i +_p b_i (-E_i)))
\end{align*}
for every $ P \in \overline{\mathcal{P}}_o ^2 $, where $\sum \nolimits ^p$ denotes the $L_p$-Minkowski sum which is taken over $E_i \in \xi _o (P)$.  If $r = 2/p-1$, then there are constants $ a \geq 0 $ and $-1 \leq \tau \leq 1 $ such that $$Z P =a \hat{\Pi} _p ^ \tau P $$
for every $ P \in \overline{\mathcal{P}}_o ^2 $.
In all other cases, $Z P = \{o\} $ for every $ P \in \overline{\mathcal{P}}_o ^2 $.
\end{lem}

Now we classify continuous, homogeneous, $SL(n)$ covariant normalized symmetric $L_p$-Blaschke valuations.
\begin{thm}\label{ncop}
Let $n \geq 3$, $p > 1$ and $p$ not an even integer. Then there exist no continuous, homogeneous, $SL(n)$ covariant normalized symmetric $L_p$-Blaschke valuations on $\mathcal{K}_o ^n$.

Let $p >1$ and $p$ not an even integer. If $ Z : \mathcal{K}_o ^2 \rightarrow \langle \mathcal{K}_c ^2,\widetilde{\#}_p \rangle$ is a continuous, homogeneous, $SL(2)$ covariant valuation, then there exists a constant $c>0$ such that
$$ZK=c \psi _{\pi /2} \widetilde{\Lambda} _c ^p K$$
for every $K \in \mathcal{K}_o ^2$.
\begin{proof}
Assume that $ Z : \mathcal{K}_o ^n \rightarrow \langle \mathcal{K}_c ^n,\widetilde{\#}_p \rangle$ is a continuous, $SL(n)$ covariant valuation which is homogeneous of degree $q$.  Lemma \ref{4.1} shows $q \leq -n+1$.

We firstly consider the cases $n \geq 3$.  If $q < -n + 1$, then $\overline{Z}_2$, introduced in Lemma \ref{coMin2}, is a $SL(n)$ contravariant $L_p$-Minkowski valuation which is homogeneous of degree $r>n-1$.  By Lemma \ref{ludwigcop}, we have $\overline{Z}_2 P = \{o\}$ for every $P \in \overline{\mathcal{P}}_o ^n$.  If $q = -n + 1$, $\overline{Z}_1$, introduced in Lemma \ref{coMin1}, is a $SL(n)$ contravariant $L_p$-Minkowski valuation which is homogeneous of degree $n-1$.  By Lemma \ref{ludwigcop}, $\overline{Z}_1 P = \{o\}$ for every $ P \in \overline{\mathcal{P}}_o ^n$.

Combined with the injectivity relation of the $L_p$-cosine transform (\ref{inj}), all cases $q \leq -n+1$ imply that
\begin{align*}
\frac{{S_p(ZP,\cdot)}}{V(ZP)} = 0
\end{align*}
for every $P \in \overline{\mathcal{P}}_o ^n$.  It is a contradiction to the existence of continuous, homogeneous, $SL(n)$ covariant normalized symmetric $L_p$-Blaschke valuations on $\mathcal{K}_o ^n$.

Next we consider the case $n=2$.  If $q < -1,~q \neq -2/p-1$, then $\overline{Z}_2$, introduced in Lemma \ref{coMin2}, is a $SL(2)$ contravariant $L_p$-Minkowski valuation which is homogeneous of degree $r>1,~r \neq 2/p+1$.  By Lemma \ref{ludwigcop}, we have $\overline{Z}_2 P = \{o\}$ for every $P \in \overline{\mathcal{P}}_o ^2$.  Combined with the injectivity relation of the $L_p$-cosine transform (\ref{inj}), we get $\frac{{S_p(ZP,\cdot)}}{V(ZP)} = 0$.  That is a contradiction.

If $q = -2/p-1$, then $\overline{Z}_2$, introduced in Lemma \ref{coMin2}, is a $SL(2)$ contravariant $L_p$-Minkowski valuation which is homogeneous of degree $2/p+1$.  By Lemma \ref{ludwigcop}, there are constants $a \geq 0$ and $-1 \leq \tau \leq 1$ such that
\begin{align*}
\overline{Z}_2 P = a \psi _{\pi /2} M_p ^ \tau P
\end{align*}
for every $P \in \overline{\mathcal{P}}_o ^2$.  Thus, $\psi _{- \pi /2} \overline{Z}_2 P = a M_p ^ \tau P$ for every $P \in \mathcal{P}_o ^2$.  Assumption that $Z$ does not contain $\{o \}$ in its range gives $a>0$.  Since $\overline{Z}_2 P$ is origin-symmetric, we get $\tau =0$.  Thus,
$\psi _{- \pi /2} (\frac{\Pi _p ZP}{V(ZP)^{1/p}}) =a M_p ^0 P$ for every $P \in \mathcal{P}_o ^2$.  Since the operators $\psi _{-\pi /2},~ \frac{\Pi _p Z}{V^{1/p}}$ and $\Gamma _p$ are continuous on $\mathcal{K}_o ^2$, and $\mathcal{P}_o ^2$ is dense in $\mathcal{K}_o ^2$, we obtain
\begin{align*}
\psi _{- \pi /2} (\frac{{\Pi _p ZK}}{V(ZK)^{1/p}}) =a M_p ^0 K
\end{align*}
for every $K \in \mathcal{K}_o ^2$.  By rewriting this in terms of the $L_p$-cosine transforms (via relation (\ref{2.11}) and $(c_{n,p}V(K))^\frac{1}{p}\Gamma _p K = M_p ^0 K$), we get
\begin{align*}
C_p \frac{{S_p(ZK,\cdot)}}{V(ZK)} (\psi _{\pi /2} x)
= b C_p (\frac{{1}}{2} \rho _K (\cdot)^{n+p} + \frac{{1}}{2} \rho _{-K} (\cdot)^{n+p})(x)
\end{align*}
for a suitable constant $b>0$.  Since
\begin{align*}
C_p \frac{{S_p(\psi _{- \pi /2}ZK,\cdot)}}{V(\psi _{- \pi /2}ZK)} (x) = C_p \frac{{S_p(ZK,\cdot)}}{V(ZK)} (\psi _{\pi /2} x)
\end{align*}
(by relation (\ref{ncp})), the injectivity property (\ref{inj}) and the definition of the normalized symmetric $L_p$-curvature image operator finally show
\begin{align*}
\psi _{- \pi /2} ZK = c \widetilde{\Lambda} _c ^p K
\end{align*}
for a suitable constant $c>0$.  Hence,
\begin{align*}
ZK = c \psi _{\pi /2}  \widetilde{\Lambda} _c ^p K
\end{align*}
for every $K \in \mathcal{K}_o ^2$

If $q = -1$, $\overline{Z}_1$, introduced in Lemma \ref{coMin1}, is a $SL(2)$ contravariant $L_p$-Minkowski valuation which is homogeneous of degree $1$.  By Lemma \ref{ludwigcop}, there are constants $ a_0,b_0 \geq 0 $ and $a_i,b_i \in \mathbb{R}$ with $a_0+a_i,~b_0+b_i \geq 0,~i=1,2$ such that
\begin{align*}
\overline{Z}_1 P =\psi _{\pi /2} (a_0P +_p b_0(-P)+_p \sum \nolimits ^p (a_i E_i +_p b_i (-E_i)))
\end{align*}
for every $ P \in \overline{\mathcal{P}}_o ^2 $, where $\sum \nolimits ^p$ denotes the $L_p$-Minkowski sum, and the sum is taken over $E_i \in \xi _o (P)$.  For $P_0 = [\pm e_1,\pm e_2]$, we have
\begin{align*}
\frac{{\Pi _p ZP_0}}{V(ZP_0)^{1/p}} = c \psi _{\pi /2}P_0,
\end{align*}
for a suitable $c \geq 0$.  Assumption that $Z$ does not contain $\{o \}$ in its range gives $c>0$.  For $p>1$, every support set of a $L_p$-projection body consists of precisely one point (Lemma \ref{pgz}).  However, $\psi _{\pi /2} P_0$ has a support set $[e_1,e_2]$ which does not consist of precisely one point. That is a contradiction.
\end{proof}
\end{thm}

Theorem \ref{th3.1} and Theorem \ref{ncop} now directly imply Theorem \ref{1.2}.

\section{$L_p$-Blaschke Valuations}
We firstly give the relationship between normalized symmetric $L_p$-Blaschke valuations and symmetric $L_p$-Blaschke valuations.
\begin{lem}\label{th5.1}
If $ Z : \mathcal {Q} \rightarrow \langle \mathcal{K}_c ^n,\# _p \rangle$ is a symmetric $L_p$-Blaschke valuation, then $ \widetilde{Z} : \mathcal {Q} \rightarrow \langle \mathcal{K}_c ^n,\widetilde{\#}_p \rangle$, defined by
\begin{align}\label{5.1}
\frac{{S_p(\widetilde{Z}K,\cdot)}}{V(\widetilde{Z}K)} = S_p(ZK)
\end{align}
for every $K \in \mathcal {Q}$, is a normalized symmetric $L_p$-Blaschke valuation.  Moreover, $ \widetilde{Z}$ is continuous if $Z$ is continuous, $ \widetilde{Z}$ is $SL(n)$ covariant (or contravariant) if $Z$ is $SL(n)$ covariant (or contravariant respectively), and $\widetilde{Z}$ is homogeneous of degree $q(p-n)/p$ if $Z$ is homogeneous of degree $q$.
\begin{proof}
Since $Z$ is a symmetric $L_p$-Blaschke valuation,
\begin{align*}
    S_p(Z(K \cup L),\cdot) + S_p(Z(K \cap L),\cdot) = S_p(ZK,\cdot) + S_p(ZL,\cdot),
\end{align*}
whenever $K,L,K \cup L,K \cap L \in \mathcal {Q}$.  By the definition of $\widetilde{Z}$ and normalized $L_p$-Blaschke sum, $\widetilde{Z}$ is a normalized symmetric $L_p$-Blaschke valuation.

We can prove continuity of $\widetilde{Z}$ in a similar way to show continuity of the normalized symmetric $L_p$-curvature image.  But because of the existence of $ZK$, we can prove it in an easier way (without using Lemma \ref{th4.2}).

By the uniqueness of the volume-normalized even $L_p$-Minkowski problem, we can rewrite relation (\ref{5.1}) as
\begin{align}\label{5.2}
\widetilde{Z}K = V(ZK)^{-1/p} ZK
\end{align}
for every $K \in \mathcal{K}^n$.  Since $V(ZK)>0$, if $ZK_i \to ZK$,
\begin{align*}
\widetilde{Z}K_i = V(ZK_i)^{-1/p} ZK_i \to V(ZK)^{-1/p} ZK = \widetilde{Z}K.
\end{align*}
Thus, if $Z$ is continuous, then $\widetilde{Z}$ is continuous.

If $Z(\lambda K) = \lambda ^q ZK$, for every $\lambda >0$, then
\begin{align*}
\widetilde{Z} (\lambda K) = V(Z\lambda K)^{-1/p} Z\lambda K = \lambda ^{q(p-n)/p} V(ZK)^{-1/p} ZK = \lambda ^{q(p-n)/p} \widetilde{Z}K.
\end{align*}
Thus, if $Z$ is homogeneous of degree $q$, $\widetilde{Z}$ is homogeneous of degree $q(p-n)/p$.

The proof of covariance or contravariance of $\widetilde{Z}$ is similar to the proof of homogeneity.
\end{proof}
\end{lem}

Lemma \ref{th5.1} introduces a map from the space of symmetric $L_p$-Blaschke valuations to the space of normalized symmetric $L_p$-Blaschke valuations, and the continuity, homogeneity or $SL(n)$ covariance (or contravariance) of symmetric $L_p$-Blaschke valuations are inherited by the corresponding normalized cases.
For $p \neq n$, the relation (\ref{5.1}) can also be rewritten as
\begin{align}\label{5.3}
V(\widetilde{Z}K)^{1/(p-n)}\widetilde{Z}K = ZK
\end{align}
for every $K \in \mathcal {Q}$.  Then we get the following lemma in a similar way.  Hence, the map is a bijection and these properties are also transferred by the inverse map.

\begin{lem}\label{th5.2}
If $ \widetilde{Z} : \mathcal {Q} \rightarrow \langle \mathcal{K}_c ^n,\widetilde{\#}_p \rangle$ is a normalized symmetric $L_p$-Blaschke valuation, $p \neq n$, then $ Z : \mathcal {Q} \rightarrow \langle \mathcal{K}_c ^n,\# _p \rangle$, defined by
\begin{align}
ZK = V(\widetilde{Z}K)^{1/(p-n)} \widetilde{Z}K
\end{align}
for every $K \in \mathcal {Q}$, is a symmetric $L_p$-Blaschke valuation.  Moreover, $Z$ is continuous if $ \widetilde{Z}$ is continuous, $Z$ is $SL(n)$ covariant (or contravariant) if $ \widetilde{Z}$ is $SL(n)$ covariant (or contravariant respectively), and $Z$ is homogeneous of degree $qp/(p-n)$ if $ \widetilde{Z}$ is homogeneous of degree $q$.
\end{lem}

Lemma \ref{th5.1}, Lemma \ref{th5.2} together with Theorem \ref{1.1} (or Theorem \ref{th3.1} as well as Theorem \ref{npbt}) provide a classification of continuous, homogeneous $SL(n)$ contravariant symmetric $L_p$-Blaschke valuations on $\mathcal{K}_o ^n$.
\begin{thm}\label{th5.3}
    For $n \geq 2$, $p > 1,~p \neq n$ and $p$ not an even integer, a map $Z : \mathcal{K}_o ^n \rightarrow \langle \mathcal{K}_c ^n,\# _p \rangle$ is a continuous, homogeneous, $SL(n)$ contravariant symmetric $L_p$-Blaschke valuation, if and only if there exists a constant $c>0$ such that
$$ZK=c \Lambda _c ^p K$$
for every $K \in \mathcal{K}_o ^n$.
\begin{proof}
Since $Z$ is a continuous, homogeneous $SL(n)$ contravariant symmetric $L_p$-Blaschke valuation, $\widetilde{Z}$ defined in Lemma \ref{th5.1} is a continuous, homogeneous $SL(n)$ contravariant normalized symmetric $L_p$-Blaschke valuation.
Theorem \ref{npbt} implies that there exists a constant $c>0$ such that
$$\widetilde{Z} K= c \widetilde{\Lambda} _c ^p K$$ for every $K \in \mathcal{K}_o ^n$.
Note that $\Lambda _c ^p K = V(\widetilde{\Lambda} _c ^pK)^{1/(p-n)} \widetilde{\Lambda} _c ^p K$. By relation (\ref{5.3}),
\begin{align}
ZK = V(\widetilde{Z}K)^{1/(p-n)} \widetilde{Z}K = V(c \widetilde{\Lambda} _c ^p K)^{1/(p-n)} c \widetilde{\Lambda} _c ^p K = c^{p/(p-n)} \Lambda _c ^p K
\end{align}
for every $K \in \mathcal{K}_o ^n$.

On the other hand, Theorem \ref{th3.1} implies that $\widetilde{\Lambda} _c ^p K$ is a continuous, homogeneous $SL(n)$ contravariant normalized symmetric $L_p$-Blaschke valuation.  Then, $\Lambda _c ^p K$ is a continuous, homogeneous, $SL(n)$ contravariant symmetric $L_p$-Blaschke valuation by Lemma \ref{th5.2}.
\end{proof}
\end{thm}

Lemma \ref{th5.1} and Lemma \ref{th5.2} together with Theorem \ref{1.2} (or Theorem \ref{th3.1} as well as Theorem \ref{ncop}) imply the following theorem.

\begin{thm}\label{th5.5}
Let $n \geq 3$, $p >1$ and $p$ not an even integer. Then there exist no continuous, homogeneous, $SL(n)$ covariant symmetric $L_p$-Blaschke valuations on $\mathcal{K}_o ^n$.

Let $p >1$ and $p$ not an even integer. If $ Z : \mathcal{K}_o ^2 \rightarrow \langle \mathcal{K}_c ^2,\# _p \rangle$ is a continuous, homogeneous, $SL(2)$ covariant symmetric $L_p$-Blaschke valuation, then there exists a constant $c>0$ such that
$$ZK=c \psi _{\pi /2} \Lambda _c ^p K$$
for every $K \in \mathcal{K}_o ^2$.
\begin{proof}
For $n \geq 3$, we argue by contradiction.  Assume that $Z$ is a continuous, homogeneous, $SL(n)$ covariant symmetric $L_p$-Blaschke valuation, $\widetilde{Z}$ defined in Lemma \ref{th5.1} is a continuous, homogeneous, $SL(n)$ covariant normalized symmetric $L_p$-Blaschke valuation.
But Theorem \ref{ncop} implies that there are no continuous, homogeneous, $SL(n)$ covariant normalized symmetric $L_p$-Blaschke valuations on $\mathcal{K}_o ^n$.  That is a contradiction.

For $n=2$, the proof is almost the same as in Theorem \ref{th5.3}.
\end{proof}
\end{thm}

\bibliographystyle{amsalpha}

\end{document}